\documentclass [11pt]{amsart}

\usepackage {amsmath, amssymb, amscd, graphicx, color, url,epsf}
\usepackage[all,knot,arc]{xy}
\usepackage[text={6.5in,9in},centering,letterpaper,dvips]{geometry}

\newtheorem {theorem}{Theorem}[section]
\newtheorem {lemma}[theorem]{Lemma}
\newtheorem {prop}[theorem]{Proposition}
\newtheorem {corollary}[theorem]{Corollary}

\newtheorem {definition}[theorem]{Definition}

\theoremstyle{remark}
\newtheorem {remark}[theorem]{Remark}

\newcommand\Z{\mathbb{Z}}

\newcommand\Ta{\mathbb{T}_\alpha}
\newcommand\Tb{\mathbb{T}_\beta}
\newcommand\Tg{\mathbb{T}_\gamma}

\newcommand{\torsion}{\textit{torsion}}

\DeclareMathOperator\rank{rank}

\DeclareMathOperator\im{image}
\DeclareMathOperator\Hom{Hom}

\def\Sym{\mathrm{Sym}}

\def\zz {{\mathbb{Z}}}

\def\F {{\mathbb{F}}}
\def\pidelta{\pi_{\Delta}}
\def\pisigma{\pi_{\Sigma}}
\def\del {{\partial}}

\def\spc {{\operatorname{Spin^c}}}
\def\rk {{\operatorname{rank}}}
\def\fin\qedhere
\def\pr {{\text{pr}}}

\def\hf{\widehat{HF}}

\def\br {\operatorname{br}}
\def\s{\mathfrak s}

\def\alphas{\boldsymbol{\alpha}}
\def\betas{\boldsymbol{\beta}}
\def\gammas{\boldsymbol{\gamma}}

\newcommand{\zees}{\mathbf{z}}
\newcommand{\wees}{\mathbf{w}}
\newcommand{\ttt}{\mathfrak{t}}
\newcommand{\Fmap}{\hat{F}}

\newcommand{\calh}{\mathcal{H}}

{\makeatletter\gdef\reallynopagebreak{\nopagebreak\@nobreaktrue}}

\begin{document}

\title{Combinatorial cobordism maps in hat Heegaard Floer theory}

\author[Robert Lipshitz]{Robert Lipshitz}
\thanks {RL was supported by an NSF Mathematical Sciences Postdoctoral Research Fellowship.}
\address {Department of Mathematics, Columbia University\\ New York, NY 10027}
\email {lipshitz@math.columbia.edu}

\author[Ciprian Manolescu]{Ciprian Manolescu}
\thanks {CM was supported by a Clay Research Fellowship.}
\address {Department of Mathematics, Columbia University\\ New York, NY 10027}
\email {cm@math.columbia.edu}

\author[Jiajun Wang]{Jiajun Wang}
\thanks{JW was partially supported by the NSF Holomorphic Curves FRG grant 0244663.}
\address{Department of Mathematics, California Institute of Technology\\ 
Pasadena, CA 91125}
\email {wjiajun@caltech.edu}

\begin {abstract} In a previous paper, Sarkar and the third author gave a 
combinatorial description of the hat version of Heegaard Floer homology 
for three-manifolds. Given a cobordism between two connected 
three-manifolds, there is an induced map between their Heegaard Floer 
homologies. Assume that the first homology group of each boundary 
component surjects onto the first homology group of the cobordism (modulo 
torsion). Under this assumption, we present a procedure for finding the 
rank of the induced Heegaard Floer map combinatorially, in the hat 
version. \end {abstract}

\maketitle

\section{Introduction}
\label{sec:Intro}

In their papers \cite{HolDisk}, \cite{HolDiskTwo}, \cite{HolDiskFour}, Ozsv\'ath and Szab\'o constructed a decorated topological quantum field theory (TQFT) in $3+1$ dimensions, called Heegaard Floer theory. (Strictly speaking, the axioms of a TQFT need to be altered slightly.) In its simplest version (called hat), to a closed, connected, oriented three-manifold $Y$ and a $\spc$ structure $\s$ on $Y$ one associates a vector space $\hf(Y,\s)$ over the field $\F=\Z/2\Z.$ Also, to a connected, oriented four-dimensional cobordism from $Y_1$ to $Y_2$ decorated with a $\spc$ structure $\ttt$, one associates a map 
$$ \Fmap_{W, \ttt}: \hf(Y_1, \ttt|_{Y_1}) \to \hf(Y_2, \ttt|_{Y_2}).$$

The maps $\Fmap_{W,\ttt}$ can be used to detect exotic smooth structures
on $4$-manifolds with boundary. For example, this can be seen by
considering the nucleus $X=N(2)$ of the elliptic surface $E(2)=K3,$
i.e. a regular neighborhood of a cusp fiber and a section,
cf. \cite[p. 74]{GS}. Let $X'= N(2)_p$ be the result of a log
transform with multiplicity $p$ ($p > 1$, odd) on a regular fiber $T^2
\subset X$, cf. \cite[Section 3.3]{GS}. Then $X$ and $X'$ are
homeomorphic 4-manifolds (with $\pi_1 = 1$), having as boundary the
Brieskorn sphere $\Sigma(2,3,11).$ However, they are not
diffeomorphic: this can be shown using the Donaldson or Seiberg-Witten
invariants (see \cite{Lisca}, \cite{MM}, \cite{FS}), but also by
comparing the hat Heegaard Floer invariants $\Fmap_{W, \ttt}$
and $\Fmap_{W',\ttt}$, where $W$ and $W'$ are the cobordisms
from $S^3$ to $\Sigma(2,3,11)$ obtained by deleting a $4$-ball from $X$
and $X'$, respectively. Indeed, the arguments of Fintushel-Stern \cite{FS} and Szab\'o-Stipsicz \cite{SS} can be easily adapted to show that $W$ and $W'$ have different hat Heegaard Floer invariants; one needs to use the computation of $\hf (\Sigma(2,3,11)),$ due to Ozsv\'ath-Szab\'o \cite[p.47]{OSabs}, and the rational blow-down formula of Roberts \cite{Roberts}. (It is worth noting that the maps $\Fmap$ give no nontrivial information for closed 4-manifolds, cf. \cite{HolDiskFour}; exotic structures on those can be detected with the mixed Heegaard Floer invariants of \cite{HolDiskFour}.)

The original definitions of the vector spaces $\hf$ and the maps
$\Fmap$ involved counting pseudoholomorphic disks and triangles
in symmetric products of Riemann surfaces; the Riemann surfaces are related to the three-manifolds and cobordisms involved via Heegaard diagrams. In \cite{SarkarWang}, Sarkar and the third author showed that every three-manifold admits a Heegaard diagram that is nice in the following sense: the curves split the diagram into elementary domains, all but one of which are bigons or rectangles. Using such a diagram, holomorphic disks in the symmetric product can be counted combinatorially, and the result is a combinatorial description of $\hf(Y)$ for any $Y,$ as well as of the hat version of Heegaard Floer homology of null homologous knots and links in any three-manifold $Y$. A similar result was obtained in \cite{MOS} for all versions of the Heegaard Floer homology of knots and links in the three-sphere.

The goal of this paper is to give a combinatorial procedure for calculating the ranks of the maps $\Fmap_{W,\s}$ when $W$ is a cobordism between $Y_1$ and $Y_2$ with the property that the induced maps $H_1(Y_1; \zz)/\torsion \to H_1(W; \zz)/\torsion$ and $H_1(Y_2; \zz)/\torsion \to H_1(W; \zz)/\torsion$ are surjective. Note that this case includes all cobordisms for which $H_1(W; \zz)$ is torsion, as well as all those consisting of only 2-handle additions.

Roughly, the computation of the ranks of $\Fmap_{W,\s}$ goes as follows. The cobordism $W$ is decomposed into a sequence of one-handle additions, two-handle additions, and three-handle additions. Using the homological hypotheses on the cobordism and the $(H_1/\torsion)$-action on the Heegaard Floer groups we reduce the problem to the case of a cobordism map corresponding to two-handle additions only. Then, given a cobordism made of two-handles, we show that it can be represented by a multi-pointed triple Heegaard diagram of a special form, in which all elementary domains that do not contain basepoints are bigons, triangles, or rectangles. In such diagrams all holomorphic triangles of Maslov index zero can be counted algorithmically, thus giving a combinatorial description of the map on $\hf.$

We remark that in order to turn $\hf$ into a fully combinatorial TQFT (at least for cobordisms satisfying our hypothesis), one ingredient is missing: naturality. Given two different nice diagrams for a three-manifold, the results of \cite{HolDisk} show that the resulting groups $\hf$ are isomorphic. However, there is not yet a combinatorial description of this isomorphism. Thus, while the results of this paper give an algorithmic procedure for computing the rank of a map $\Fmap_{W,\s},$ the map itself is determined combinatorially only up to automorphisms of the image and the target. In fact, if one were to establish naturality, then one could automatically remove the assumption on the maps on $H_1/\torsion$, and compute $\Fmap_{W,\s}$ for any $W$, simply by composing the maps induced by the two-handle additions (computed in this paper) with the ones induced by the one- and three-handle additions, which are combinatorial by definition, cf. \cite{HolDiskFour}.

The paper is organized as follows. In Section~\ref{sec:Triangles}, we define a multi-pointed triple Heegaard diagram to be nice if all non-punctured elementary domains are bigons, triangles, or rectangles, and show that in a nice diagram holomorphic triangles can be counted combinatorially.\footnote{Sucharit Sarkar has independently obtained this result (Proposition~\ref{EasternOrthodox} below) in \cite{SarkarIndex}, using slightly different methods.} We then turn to the description of the map induced by two-handle additions. For the sake of clarity, in Section~\ref{sec:Two} we explain in detail the case of adding a single two-handle: we show that its addition can be represented by a nice triple Heegaard diagram with a single basepoint and, therefore, the induced map on $\hf$ admits a combinatorial description. We then explain how to modify the arguments to work in the case of several two-handle additions.  This modification uses triple Heegaard diagrams with several basepoints. In Section~\ref{sec:Last}, we discuss the additions of one- and three-handles, and put the various steps together. Finally, in Section~\ref{Sec:Example} we present the example of $+1$ surgery on the trefoil.

Throughout the paper all homology groups are taken with coefficients in $\F = \zz/2\zz,$ unless otherwise noted.

\subsection*{Acknowledgments}

We would like to thank Peter Ozsv\'ath and Zolt\'an Szab\'o for helpful conversations and encouragement. In particular, several key ideas in the proof were suggested to us by Peter Ozsv\'ath.

This work was done while the third author was an exchange graduate student at Columbia University. He is grateful to the Columbia math department for its hospitality. He would also like to thank his advisors, Robion Kirby and Peter Ozsv{\'a}th, for their continuous guidance and support.

Finally, we would like to thank the referees for many helpful comments, and particularly for finding a critical error in Section \ref{sec:Last} of a previous version of this paper.

\newpage
\section{Holomorphic triangles in nice triple Heegaard diagrams}
\label{sec:Triangles}

The goal of this section is to show that under an appropriate condition (``niceness'') on triple Heegaard diagrams, the counts of holomorphic triangles in the symmetric product are combinatorial.

\subsection{Preliminaries} \label{sec:prels}
We start by reviewing some facts from Heegaard Floer theory. A triple Heegaard diagram $\calh = (\Sigma, \alphas, \betas, \gammas)$ consists of a surface $\Sigma$ of genus $g$ together with three $(g+k)$-tuples of pairwise disjoint embedded curves $\alphas = \{\alpha_1, \dots, \alpha_{g+k}\},$ $ \betas=\{\beta_1, \dots, \beta_{g+k}\}, \gammas= \{\gamma_1, \dots, \gamma_{g+k}\}$ in $\Sigma$ such that the span of each $(g+k)$-tuple of curves in $H_1(\Sigma)$ is $g$-dimensional. If we forget one set of curves (for example $\gammas$), the result is an (ordinary) Heegaard diagram $(\Sigma, \alphas, \betas).$

By the condition on the spans, $\Sigma\setminus\alphas$, $\Sigma\setminus\betas$ and $\Sigma\setminus\gammas$ each has $k+1$ connected components. By a multi-pointed triple Heegaard diagram $(\calh,\zees)$, then, we mean a triple Heegaard diagram $\calh$ as above together with a set $\zees=\{z_1,\ldots,z_{k+1}\} \subset \Sigma$ of $k+1$ points in $\Sigma$ so that exactly one $z_i$ lies in each connected component of $\Sigma\setminus\alphas$, $\Sigma\setminus\betas$ and $\Sigma\setminus\gammas$.

To a Heegaard diagram $(\Sigma, \alphas,\betas)$ one can associate a
three-manifold $Y_{\alpha,\beta}$. To a triple Heegaard diagram
$(\Sigma,\alphas,\betas,\gammas)$, in addition to the three-manifolds
$Y_{\alpha,\beta}$, $Y_{\beta,\gamma}$ and $Y_{\alpha,\gamma}$, one
can associate a four-manifold $W_{\alpha,\beta,\gamma}$ such that
$\partial W_{\alpha,\beta,\gamma}=-Y_{\alpha,\beta}\cup
-Y_{\beta,\gamma} \cup Y_{\alpha,\gamma}$; see~\cite{HolDisk}.

Associated to a three-manifold $Y$ is the Heegaard Floer homology
group $\hf(Y)$. This was defined using a Heegaard diagram with a
single basepoint in~\cite{HolDisk}. In~\cite{Links}, Ozsv\'ath and
Szab\'o associated to the data $(\Sigma,\alphas,\betas,\zees)$, called
a multi-pointed Heegaard diagram, a Floer homology group
$\hf(\Sigma,\alphas,\betas,\zees)$ by counting holomorphic disks in
$\Sym^{g+k}(\Sigma\setminus\zees)$ with boundary on the tori
$\Ta=\alpha_1\times\cdots\times\alpha_{g+k}$ and
$\Tb=\beta_1\times\cdots\times\beta_{g+k}.$ It is not hard to show
that
\begin{equation}\label{ManyPtsIso} 
\hf(\Sigma,\alphas,\betas,\zees) \cong \hf\left(Y_{\alpha,\beta}\right)\otimes H_*(T^k). 
\end{equation} 
(Here, $T^k$ is the $k$-torus, and $H_*(T^k)$ means ordinary
(singular) homology.) The decomposition \eqref{ManyPtsIso} is not canonical: it depends
on a choice of paths in $\Sigma$ connecting $z_i$ to $z_1$ for
$i=2,\cdots,k+1$.

The Heegaard Floer homology groups decompose as a direct sum over $\spc$-structures on $Y_{\alpha,\beta}$, 
\[ 
\hf(Y_{\alpha,\beta})\cong 
\bigoplus_{\s\in\spc(Y_{\alpha,\beta})}\hf(Y_{\alpha,\beta},\s). 
\] 
More generally, there is a decomposition 
\[ 
\hf(\Sigma,\alphas,\betas,\zees)\cong 
\bigoplus_{\s\in\spc(Y_{\alpha,\beta})}\hf(\Sigma,\alphas,\betas,\zees,\s)\cong 
\bigoplus_{\s\in\spc(Y_{\alpha,\beta})}\left(\hf(Y_{\alpha,\beta}, 
\s)\otimes H_*(T^k)\right). 
\]

Associated to the triple Heegaard diagram
$(\Sigma,\alphas,\betas,\gammas,\zees)$ together with a
$\spc$-structure $\ttt$ on $W_{\alpha,\beta,\gamma}$, is a map
\begin{equation}\label{eq:FMaps}
\Fmap_{\Sigma,\alphas,\betas,\gammas,\zees,\ttt}:
\hf(\Sigma,\alphas,\betas,\zees,\ttt|_{Y_{\alpha,\beta}} ) 
\otimes\hf(\Sigma,\betas,\gammas,\zees, \ttt|_{Y_{\beta,\gamma}})\to
\hf(\Sigma,\alphas,\gammas,\zees,\ttt|_{Y_{\alpha,\gamma}}).
\end{equation}
The definition involves counting holomorphic triangles in
$\Sym^{g+k}(\Sigma\setminus \zees)$ with boundary on $\Ta$, $\Tb$ and
$\Tg$, cf.~\cite{HolDisk} and~\cite{Links}.

Two triple Heegaard diagrams are called \emph{strongly equivalent} if
they differ by a sequence of isotopies and handleslides. It follows from \cite[Proposition 8.14]{HolDisk}, the associativity theorem \cite[Theorem 8.16]{HolDisk}, and the definition of the handleslide isomorphisms that
strongly equivalent triple Heegaard diagrams induce the same map on homology.

 Call a
$(k+1)$-pointed triple Heegaard diagram \emph{split} if it is obtained
from a singly-pointed Heegaard triple diagram
$(\Sigma',\alphas',\betas',\gammas',z')$ by attaching (by connect sum)
$k$ spheres with one basepoint and three isotopic curves (one
alpha, one beta and one gamma) each, to the component of
$\Sigma'\setminus(\alphas'\cup\betas'\cup\gammas')$ containing $z'$. We
call $(\Sigma',\alphas',\betas',\gammas',z')$ the \emph{reduction} of
the split diagram. 

The following lemma is a variant of~\cite[Proposition 3.3]{Links}. 
\begin{lemma}\label{lemma:splitDiagrams}
Every triple Heegaard diagram $(\Sigma,\alphas,\betas,\gammas,\zees)$
is strongly equivalent to a split one.
\end{lemma}
\begin{proof}
  Reorder the alpha circles so that $\alpha_1,\dots,\alpha_g$ are
  linearly independent. Let $R$ be the connected component of
  $\Sigma\setminus\alphas$ containing $z_2$. Since
  $\alpha_1,\dots,\alpha_g$ are linearly independent, one of the
  curves $\alpha_i$, $i>g$, must appear in the boundary of $R$ with
  multiplicity exactly $1$. By handlesliding this curve over the other
  boundary components of $R$, we can arrange that the resulting
  $\alpha_i$ bounds a disk containing only $z_2$. Repeat this process
  with the other $z_i$, $i=3,\dots,k+1$, being sure to use a different
  alpha curve in the role of $\alpha_i$ at each step. We reorder the
  curves $\alpha_{g+1},\dots,\alpha_{g+k}$ so that $\alpha_i$
  encircles $z_{i-g+1}$. Now repeat the entire process for the beta
  and gamma curves.  Finally, choose a path $\zeta_i$ in $\Sigma$ from
  $z_i$ to $z_1$ for each $i=2,\dots,k+1$. Move the configuration
  $(z_i,\alpha_i,\beta_i,\gamma_i)$ along the path $\zeta_i$ by
  handlesliding (around $\alpha_i$, $\beta_i$ or $\gamma_i$) the other
  alpha, beta and gamma curves that are encountered along the path
  $\zeta_i$. The result is a split triple Heegaard diagram. Note that
  its reduction is obtained from the original diagram
  $(\Sigma,\alphas,\betas,\gammas,\zees)$ by simply forgetting some of
  the curves and basepoints.
\end{proof}

The maps from~\eqref{eq:FMaps} are compatible with the
isomorphism~\eqref{ManyPtsIso}, in the following sense.  Given a
triple Heegaard diagram $(\Sigma,\alphas,\betas,\gammas,\zees)$, let
$(\Sigma',\alphas',\betas',\gammas',z')$ be the reduction of a split
diagram strongly equivalent to
$(\Sigma,\alphas,\betas,\gammas,\zees)$.  Then the following diagram
commutes:
\begin{equation}\label{diagram:commutes}
\xymatrix{
  \left(\hf\left(Y_{\alpha,\beta},\ttt|_{Y_{\alpha,\beta}}\right)\otimes
    H_*\left(T^k\right)\right) \otimes
  \left(\hf\left(Y_{\beta,\gamma},\ttt|_{Y_{\beta,\gamma}}\right)\otimes
    H_*\left(T^k\right)\right) \ar[r]^{} \ar[d]^{\cong} &
  \hf\left(Y_{\alpha,\gamma}, \ttt|_{Y_{\alpha,\gamma}}\right)\otimes
  H_*\left(T^k\right)\ar[d]^{\cong}\\
  \hf(\Sigma,\alphas,\betas,\zees,\ttt|_{Y_{\alpha,\beta}} )
  \otimes\hf(\Sigma,\betas,\gammas,\zees, \ttt|_{Y_{\beta,\gamma}})
  \ar[r]^(.62){\Fmap_{\Sigma,\alphas,\betas,\gammas,\zees,\ttt}} &
  \hf(\Sigma,\alphas,\gammas,\zees,\ttt|_{Y_{\alpha,\gamma}}).  } 
\end{equation}
Here, the map in the first row is
$\Fmap_{\Sigma',\alphas',\betas',\gammas',z',\ttt}$ on the
$\hf$-factors and the usual intersection product
$H_*(T^k)\otimes H_*(T^k)\to H_*(T^k)$ on the $H_*(T^k)$-factors. The
vertical isomorphisms are induced by the strong equivalence.  The
proof that the diagram commutes follows from the same ideas as in
\cite{Links}: Each of the $k$ spherical pieces in the split diagram
contributes a $H_*(S^1)$ to the $H_*(T^k) =
\bigl(H_*(S^1)\bigr)^{\otimes k}$ factors above; moreover, a local
computation shows that the triangles induce intersection product maps
$H_*(S^1) \otimes H_*(S^1) \to H_*(S^1),$ which tensored together give
the intersection product on $H_*(T^k).$

\begin{remark}It is tempting to assert that the map
$\Fmap_{\Sigma,\alphas,\betas,\gammas,z,\ttt}$ induced by a
singly-pointed triple Heegaard diagram depends only on
$W_{\alpha,\beta,\gamma}$ and $\ttt$. However, this seems not to be known.
\end{remark}

\subsection{Index formulas} 

Fix a triple Heegaard diagram $\mathcal{H}=(\Sigma,\alphas,\betas,\gammas)$ as above. The complement of the $3(g+k)$ curves in $\Sigma$ has several connected components, which we denote by $D_1, \dots, D_N$ and call {\em elementary domains}.

The {\em Euler measure} of an elementary domain $D \subset \Sigma$ is 
$$ e(D) = \chi(D) - \frac{\# \text{ vertices of } D}{4}.$$

A {\em domain} in $\Sigma$ is a two-chain $D = \sum a_iD_i$ with $a_i \in \Z.$ Its Euler measure is simply 
$$e(D) = \sum_{i=0}^N a_i e(D_i).$$

As mentioned above, the maps $\Fmap_{\Sigma,\alphas,\betas,\gammas,\zees,\ttt}$ induced by the triple Heegaard diagram $\mathcal{H}$ are defined by counting holomorphic triangles in $\Sym^{g+k}(\Sigma)$ with respect to a suitable almost complex structure. According to the cylindrical formulation from \cite{Lipshitz}, this is equivalent to counting certain holomorphic embeddings $u: S \to \Delta \times \Sigma,$ where $S$ is a Riemann surface (henceforth called the source) with some marked points on the boundary (which we call corners), and $\Delta$ is a fixed disk with three marked points on the boundary. The maps $u$ are required to satisfy certain boundary conditions, and to be generically $(g+k)$-to-$1$ when post-composed with the projection $\pidelta: \Delta \times \Sigma \to \Delta.$ More generally, we will consider such holomorphic maps $u: S \to \Delta \times \Sigma$ which are generically $m$-to-$1$ when post-composed with $\pidelta;$ these correspond to holomorphic triangles in $\Sym^m(\Sigma),$ where $m$ can be any positive integer. We will be interested in the discussion of the index from \cite{Lipshitz}.  Although this discussion was carried-out in the case $k=0$, $m=g$, it applies equally well in the case of arbitrary $k$ and $m$ with only notational changes.

In the cylindrical formulation, one works with an almost complex structure on $\Delta\times\Sigma$ so that the projection $\pidelta$ is holomorphic, and the fibers of $\pisigma$ are holomorphic.  It follows that for $u:S\to\Delta\times\Sigma$ holomorphic, $\pidelta\circ u$ is a holomorphic branched cover. The map $\pisigma\circ u$ need not be holomorphic, but since the fibers are holomorphic, $\pisigma\circ u$ is a branched map. Fix a model for $\Delta$ in which the three marked points are 
$90^{\circ}$ corners, and a conformal structure on $\Sigma$ with respect to which the intersections between alpha, beta and gamma curves are all right angles.  Since $u$ is holomorphic, the conformal structure on $S$ is induced via $\pidelta\circ u$ from the conformal structure on $\Delta$. It makes sense, therefore, to talk about branch points of $\pisigma\circ u$ on the boundary and at the corners, as well as in the interior.  Generically, while there may be branch points of $\pisigma\circ u$ on the boundary of $S$, there will not be branch points at the corners.

Suppose $u: S \to \Delta \times \Sigma$ is as above. Denote by $\pisigma:\Delta \times \Sigma \to \Sigma$ the projection to 
$\Sigma.$ There is an associated domain $D(u)$ in $\Sigma,$ where the coefficient of $D_i$ in $D(u)$ is the local multiplicity of $\pisigma \circ u$ at any point in $D_i.$  By \cite[p. 1018]{Lipshitz}, the {\em index} of the linearized $\bar \partial$ operator at the holomorphic map $u$ is given by 
\begin {equation}\label {eq:index1}
\mu(u) = 2e(D(u)) - \chi(S) + \frac{m}{2}.
\end {equation}
For simplicity, we call this the index of $u.$

Note that, by the Riemann-Hurwitz formula:
\begin {equation}\label {eq:branch} 
\chi(S) = e(D(u)) + \frac{3m}{4} - \br(u),
\end {equation}
where $\br(u)$ is the ramification index (number of branch points counted with multiplicity) of $\pisigma \circ u.$ (Here, branch points along the boundary count as half an interior branch point.) From \eqref{eq:index1} and \eqref{eq:branch} we get an alternate formula for the index:
\begin {equation}\label {eq:index2}
\mu(u) = e(D(u)) + \br(u) - \frac{m}{4}.
\end {equation}

Note that it is not obvious how to compute $br(u)$ from
$D(u)$. A combinatorial formula for the index, purely in terms of $D(u)$, was found by Sarkar in \cite{SarkarIndex}. However, we will not use it here.

\subsection{Nice triple diagrams} 
\label{sec:NiceTriDiagrams}

Fix a multi-pointed triple Heegaard diagram $(\calh,\zees)=(\Sigma,\alphas,\betas,\gammas,\zees)$. Recall that a domain is a linear combination of connected components of $\Sigma\setminus(\alphas\cup\betas\cup\gammas)$. The \emph{support} of a 
domain is the union of those components with nonzero coefficients. If the support of a domain $D$ contains at least one $z_i$ then $D$ is called {\em punctured}; otherwise it is called {\em unpunctured}.

\begin{definition}\label{def:nice}
An elementary domain is called {\bf good} if it is a bigon, a triangle, or a rectangle, and {\bf bad} otherwise. The multi-pointed triple Heegaard diagram $(\calh,\zees)$ is called {\bf nice} if every unpunctured elementary domain is good. 
\end {definition}

This is parallel to the definition of nice Heegaard diagrams (with just two sets of curves) from \cite{SarkarWang}. A multi-pointed Heegaard diagram $(\Sigma, \alphas, \betas, \zees)$ is called {\em nice} if, among the connected components of $\Sigma \setminus (\alphas \cup \betas),$ all unpunctured ones are either bigons or squares.

Note that a bigon, a triangle, and a rectangle have Euler measure $\frac{1}{2}, \frac{1}{4},$ and $0,$ respectively. Since $e$ is additive, every unpunctured positive domain (not necessarily elementary) in a nice diagram must have nonnegative Euler measure. A quick consequence of this is the following:

\begin {lemma}\label{lemma:forget}
If $(\Sigma,\alphas,\betas,\gammas,\zees)$ is a nice triple Heegaard diagram, then if we forget one set of curves (for example, $\gammas$), the resulting Heegaard diagram $(\Sigma, \alphas, \betas, \zees)$ is also nice. 
\end {lemma}

In order to define the triangle maps it is necessary to assume the triple Heegaard diagram is weakly admissible in the sense 
of~\cite[Definition 8.8]{HolDisk}. In fact, nice diagrams are automatically weakly admissible, cf. Corollary~\ref{lemma:admissible} 
below.

Our goal is to give a combinatorial description of the holomorphic triangle counts for nice triple diagrams.

\begin {prop}\label{EasternOrthodox} 
Let $(\mathcal{H},\zees)$ be a nice multi-pointed triple Heegaard diagram. Fix a generic almost complex structure $J$ on $\Delta\times\Sigma$ as in \cite[Section 10.2]{Lipshitz}. Let $u:S \to \Delta \times \Sigma$ be a $J$-holomorphic map of the kind occurring in the definition of $\Fmap_{\alphas,\betas,\gammas,\zees,\ttt}$. In particular, assume $u$ is an embedding, of index zero, and such that the image of $\pisigma \circ u$ is an unpunctured domain. Then $S$ is a disjoint union of $m$ triangles, and the restriction of $\pisigma \circ u$ to each component of $S$ is an embedding. 
\end {prop}

\begin {proof} 
Since the image of  $\pisigma \circ u$ is unpunctured and positive, we have $e(D(u)) \geq 0.$ By \eqref{eq:index1}, we get 
$$\chi(S) \geq \frac{m}{2} > 0.$$ 

This means that at least one component of $S$ is topologically a disk. Let $S_0$ be such a component. It is a polygon with $3l$ vertices. We will show that $l=1,$ and that $\pisigma \circ u|_{S_0}$ is an embedding. 

Let us first show that $S_0$ is a triangle. The index of the $\bar\del$ operator at a disconnected curve is the sum of the indices of its restrictions to each connected component. Therefore, in order for an index zero holomorphic curve to exist generically, the indices at every connected component, and in particular at $S_0,$ must be zero. Applying \eqref{eq:index1} to $u|_{S_0}$ we get 
$$ l = 2 - 4e(D(u|_{S_0})) \leq 2.$$

\begin{figure}
\center{\includegraphics[width=120pt]{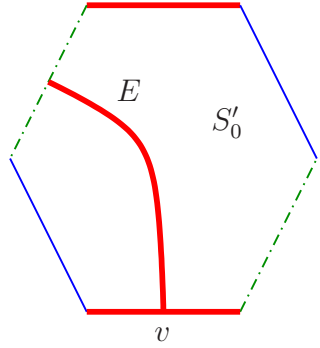}}
\caption{{\bf A hexagon component of the source.} The preimages of the alpha, beta and gamma curves (here shown as thick, thin, and interrupted lines) give an embedded graph in the source $S_0.$ A boundary branch point in the image corresponds to a valence three vertex $v$ in the source.} \label{Figure:Hexagon} 
\end{figure} 

If $l=2,$ then by \eqref{eq:index2} we have $\br(u) \leq \frac{1}{2}.$ Hence, the map
$\pisigma \circ u$ has no interior branch points.  If $\br(u)=0$ then $S_0$ is mapped
locally diffeomorphically by $\pisigma\circ u$ to $\Sigma$. The image
must have negative Euler measure, which is a contradiction. So,
suppose $\br(u)=1/2$. The preimages of the alpha, beta, and gamma
curves cut $S_0$ into several connected components. Without loss of
generality, assume that the boundary branch point is mapped to an
alpha circle. Then, along the corresponding edge of $S_0$ there is
a valence three vertex $v$, as shown in
Figure~\ref{Figure:Hexagon}. Let $E$ denote the edge in the interior
of $S_0$ meeting $v$. Since there is only one boundary branch point,
the other intersection point of the edge $E$ with $\partial S_0$ is
along the preimage of a beta or gamma circle. It follows that
one of the connected components $S_0'$ of $S_0\setminus E$ is a
hexagon or heptagon. Smoothing the vertex $v$ of $S_0'$ we obtain a
pentagon or hexagon which is mapped locally diffeomorphically by $\pisigma\circ u$ to $\Sigma$. The image, then, has negative Euler measure, again a contradiction.

Therefore, $l=1,$ so $S_0$ is a triangle. Furthermore, by \eqref{eq:index2}, $\br(u) \leq 1/4,$ which means that there are no (interior or boundary) branch points at all. Thus, just as in the hypothetical hexagon case above, the preimages of the alpha, beta, and gamma curves must cut $S_0$ into 2-gons, 3-gons, and 4-gons, all of which have nonnegative Euler measure. Since the Euler measure of $S_0$ is $1/4,$ there can be no bigons; in fact, $S_0$ must be cut into several rectangles and exactly one triangle. It is easy to see that the only possible tiling of $S_0$ of this type is as in Figure~\ref{fig:Tiling}, with several parallel preimages of segments on the alpha curves, several parallel beta segments, and several parallel gamma segments. We call the type of a segment ($\alpha$, $\beta$ or $\gamma$) its color.

\begin{figure}
\center{\includegraphics[width=150pt]{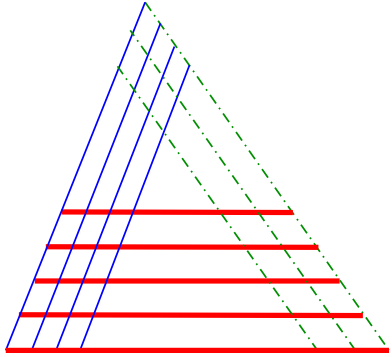}}
\caption {{\bf An embedded triangle.} The thick (red) lines are $\alpha$'s, the thin (blue) ones $\beta$'s, and the interrupted (green) ones $\gamma$'s. We first show that this is the picture in the source $S_0,$ and then that the same is true for the image of $S_0$ in $\Sigma,$ i.e. the images of all the pieces in the tiling are disjoint.}
\label{fig:Tiling}
\end{figure}

The tiling consists of one triangle and six different types of rectangles, according to the coloring of their edges in clockwise order (namely, $\alpha\beta\alpha\gamma$, $\gamma\alpha\gamma\beta$, $\beta\alpha\beta\gamma$, $\alpha\beta\alpha\beta$, $\beta\gamma\beta\gamma$, and $\gamma\alpha\gamma\alpha$). We claim that the images of the interiors of each of these rectangles by $\pisigma \circ u$ are disjoint.

Because of the coloring scheme, only rectangles of the same type can have the same image. Suppose that two different $\alpha\beta\alpha\gamma$ rectangles from $S_0$ have the same image in $\Sigma.$ (The cases $\gamma\alpha\gamma\beta$, $\beta\alpha\beta\gamma$ are exactly analogous.) Let $r_1$ and $r_2$ be the two rectangles; suppose that $r_1$ is closer to the central triangle than $r_2,$ and $r_2$ is closer to the $\alpha$ boundary of $S_0.$ Because of the way the rectangles are colored, the upper edge of $r_1$ must have the same image as the upper edge of $r_2.$ Hence the $\alpha\beta\alpha\gamma$ rectangle right above $r_1$ has the same image as the one right above $r_2.$ Iterating this argument, at some point we get that the central triangle has the same image as some $\alpha\beta\alpha\gamma$ rectangle, which is impossible.

Now suppose that two different $\alpha\beta\alpha\beta$ rectangles, $r_1$ and $r_2,$ have the same image. (The cases $\beta\gamma\beta\gamma$ and $\gamma\alpha\gamma\alpha$ are exactly analogous.) There are two cases, according to whether the upper edge of $r_1$ has the same image as the upper edge of $r_2,$ or as the lower edge of $r_2.$

Suppose first that the upper edge of $r_1$ has the same image as the upper edge of $r_2$.  By the \emph{$\beta$-height} of $r_i$ we mean the minimal number of beta arcs that an arc in $S\setminus\alphas$ starting in $r_i$, going right, and ending at a gamma arc must cross.  (The diagram is positioned in the plane as in Figure~\ref{fig:Tiling}.) Since $\pisigma\circ u$ is a local homeomorphism, and $r_1$ and $r_2$ have the same image, it is clear that the $\beta$-height of $r_1$ and the $\beta$-height of $r_2$ are equal.  By the \emph{$\alpha$-height} of $r_i$ we mean the minimal number of alpha arcs that an arc in $S\setminus\betas$ starting in $r_i,$ going up, and ending at a gamma arc must cross.  Again, it is clear that the $\alpha$-heights of $r_1$ and $r_2$ must be equal.  But this implies that $r_1$ and $r_2$ are equal.

Now, suppose that the upper edge of $r_1$ has the same image as the lower edge of $r_2$. There is a unique rectangle $R$ in $S$ with boundary contained in $\alphas\cup\betas$, containing $r_1$ and $r_2$, and with one corner the same as a corner of $r_1$ and the opposite corner the same as a corner of $r_2$. It is easy to see that $\pisigma\circ u$ maps antipodal points on the boundary of $R$ to the same point in $\Sigma$. It follows that $\pisigma\circ u|_{\partial R}$ is a two-fold covering map.  But then $\pisigma\circ u$ must have a branch point somewhere inside $R$ -- a contradiction.  See Figure~\ref{fig:Grid}.

\begin{figure}
\center{\includegraphics[width=120pt]{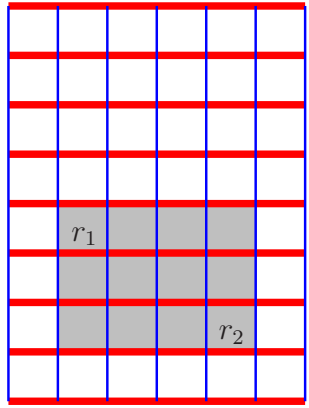}}
\caption {{\bf Two $\alpha\beta\alpha\beta$ rectangles in a grid.} If the rectangles $r_1$ and $r_2$ had the same image in $\Sigma$ with a $180^{\circ}$ turn, then the whole shaded rectangle would be mapped to $\Sigma$ with a branch point.}
\label{fig:Grid}
\end{figure}

Finally, suppose some arc $A$ on $\partial S$ has the same image as some other arc $A'$ in $S$.  If $A'$ is in the interior of $S$ then any rectangle (or triangle) adjacent to $A$ has the same image as some rectangle (or triangle) adjacent to $A'$.  We have already ruled this out.  If $A'$ is on $\partial S$, then either any rectangle (or triangle) adjacent to $A$ has the same image as some rectangle (or triangle) adjacent to $A'$ or there is a branch point somewhere on $\partial S$.  We have already ruled out both of these cases.

We have thus established that $S_0$ is an embedded triangle. By forgetting $S_0$, we obtain a holomorphic map to $\Delta \times \Sigma,$ still of index zero, but such that its post-composition with $\pidelta$ is generically $(m-1)$-to-$1$ rather than $m$-to-$1.$ The result then follows by induction on $m.$
\end {proof}

\begin{figure}
\center{\includegraphics[width=200pt]{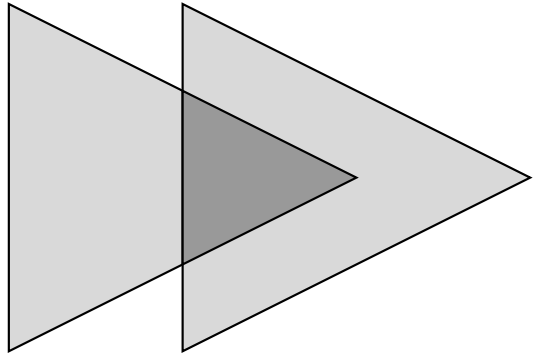}}
\caption {{\bf Head to tail overlap.} A pair of embedded triangles $T$ and $T'$ overlap ``head to tail'' if their intersection $T\cap T'$ consists of a single connected component, itself a triangle, which contains one vertex of one of the two triangles and no vertices of the other. In the picture, the two triangles are shaded; the intersection is darkly shaded. We have not colored the figure to indicate that any of the three possible coloring schemes is allowed; however, in all three cases, parallel segments in the figure are of the same type ($\alpha$, $\beta$ or $\gamma$).}
\label{fig:overlap}
\end{figure}

\begin{figure}
\center{\includegraphics[width=400pt]{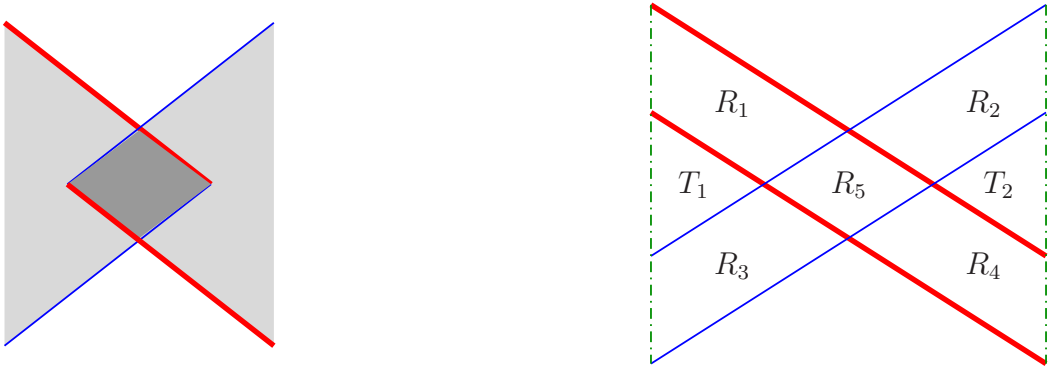}}
\caption {{\bf Head to head overlap.} If $T$ and $T'$ intersect around their $\alpha\beta$ vertices as on the left of the figure (where the intersection is shown by a darker shading), their union $T \cup T'$ does not have index zero. This can be seen using an alternate decomposition of the domain $T \cup T'$ as $T_1 \cup T_2 \cup R \cup R',$ where $R = R_1 \cup R_4 \cup R_5$ and $R' = R_2 \cup R_3 \cup R_5,$ where $T_i (i=1,2)$ and $R_j (j=1, \dots, 5)$ are the domains shown on the right. (Note that on the right, there might be overlaps in other parts of the diagram; for example, we can have the situation in Figure~\ref{fig:overlap_twice}.) The cases of $\alpha\gamma$ or $\beta\gamma$ head to head overlaps are similar.} 
\label{fig:headhead} 
\end{figure}

\begin{figure}
\center{\includegraphics[width=250pt]{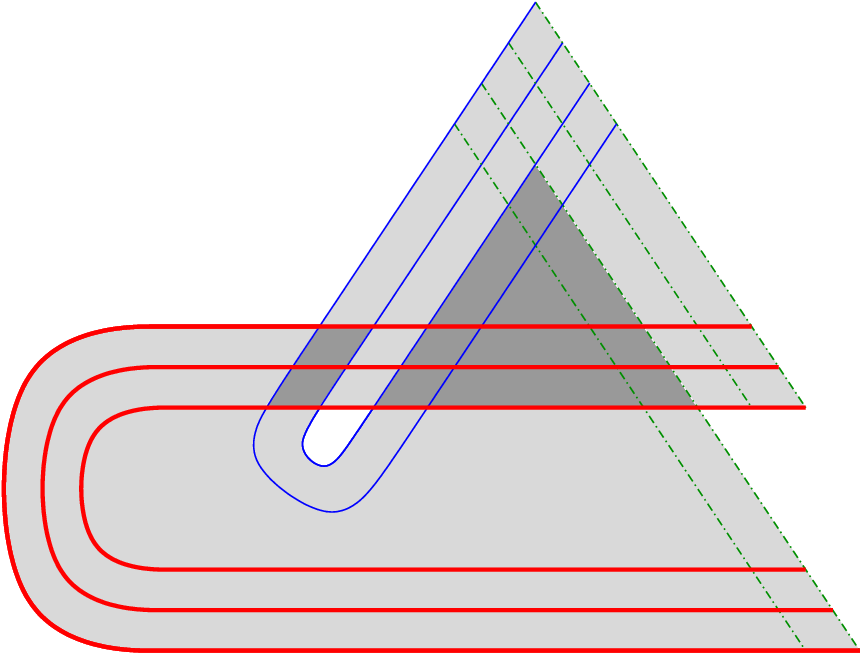}}
\caption {{\bf A double overlap.} The two triangles (each shown lightly shaded) have a (darkly shaded) overlap with two connected components.}
\label{fig:overlap_twice}
\end{figure}

Observe that, in Proposition~\ref{EasternOrthodox} above, even though each of the $m$ triangles is embedded, some of their domains may overlap. It turns out that they may do so only in a specific way, however:

\begin{lemma}\label{allow}
Suppose $A$ is an index zero homology class represented by a union of embedded holomorphic triangles, in a nice triple diagram. Suppose the union of triangles corresponds to an embedded holomorphic curve in $\Delta\times \Sigma$. Then any two triangles in $A$ are either disjoint in $\Sigma$ or overlap in $\Sigma$ ``head to tail'' as shown in Figure~\ref{fig:overlap}.
\end{lemma} 

\begin{proof}

Let $T$ and $T'$ be two of the triangles in the domain $A$. For a generic representative of $A$ to exist, the pair must also have index zero, and be embedded in $\Delta\times\Sigma$.

We already know that $T$ and $T'$ are tiled as in Figure~\ref{fig:Tiling}. This strongly restricts how $T$ and $T'$ can overlap. 

One way for $T$ and $T'$ to overlap is for $T$  to be entirely contained inside $T'$. In this case, it is not hard to see that 
the two holomorphic triangles in $\Delta\times\Sigma$ intersect in 
one interior point. Indeed, the intersection number of two holomorphic curves in a $4$-manifold is invariant in families. If we deform the Heegaard diagram so that the boundary of $T$ in $\Sigma$ is a single point (i.e., the alpha, beta and gamma circles involved in $\partial T$ intersect in an asterisk, with vertex the ``triangle'' $T$) then obviously $(T\cap T')\subset (\Delta\times\Sigma)$ is a single point. It follows that the same is true for the original triangles $T$ and $T'$.

Another way that $T$ and $T'$ might overlap is ``head to head'' as shown on the left side of Figure~\ref{fig:headhead}. It is then possible to decompose $T \cup T'$ into a pair of rectangles $R$ and $R'$, and two new embedded triangles $T$ and $T'$, as shown in Figure~\ref{fig:headhead}. An immersed rectangle in $\Sigma$ has index at least $1$, since it admits a generic holomorphic representative. So, each of $R_1$ and $R_2$ has index at least $1$. Similarly, the pair of triangles $T_1\cup T_2$ has index at least $0$.  So, by additivity of the index, the whole domain has index at least $2$ -- a contradiction.

Using these two observations, and the rulings of $T$ and $T'$, it is then elementary to check that the only possible overlap in index zero is ``head to tail'' as in Figure~\ref{fig:overlap}.
\end{proof}

It follows that, for a nice Heegaard diagram, we can combinatorially describe the generic holomorphic curves of index $0$.  If $D$ is the domain of a generic holomorphic curve of index $0$ then $\partial D$ has $m$ components, each of which bounds an embedded triangle in $\Sigma$. Each pair of triangles must either be disjoint or overlap as shown in Figure~\ref{fig:overlap}.  Any such $D$ clearly has a unique holomorphic representative with respect to a split complex structure $j_\Delta\times j_\Sigma$ on $\Delta\times\Sigma$.  Further, it is well known that these holomorphic curves are transversally cut out, and so persist if one takes a small perturbation of $j_\Delta\times j_\Sigma$. In summary, to count index zero holomorphic curves in $\Delta\times\Sigma$ with respect to a generic perturbation of the split complex structure, it suffices to count domains $D$ which are sums of $m$ embedded triangles in $\Sigma$, overlapping as allowed in the statement of Lemma~\ref{allow}.

\newpage
\section{Nice diagrams for two-handle additions}
\label{sec:Two}

In \cite[Definition 4.2]{HolDiskFour}, Ozsv\'ath and Szab\'o associate
to a four-dimensional cobordism $W$ consisting of two-handle additions
certain kinds of triple Heegaard diagrams. The cobordism $W$ from
$Y_1$ to $Y_2$ corresponds to surgery on some framed link $L \subset Y_1.$
Denote by $l$ the number of components of $L$.  Fix a basepoint in
$Y_1$. Let $B(L)$ be the union of $L$ with a path from each component
to the basepoint. The boundary of a regular neighborhood of $B(L)$ is
a genus $l$ surface, which has a subset identified with $l$
punctured tori $F_i$, one for each link component. A singly-pointed
triple Heegaard diagram $(\Sigma,\alphas,\betas,\gammas,z)$ is called
\emph{subordinate to $B(L)$} if
\begin{itemize}
\item
  $(\Sigma,\{\alpha_1,\dots,\alpha_g\},\{\beta_1,\dots,\beta_{g-l}\})$
  describes the complement of $B(L)$,
\item $\gamma_1,\dots,\gamma_{g-l}$ are small isotopic translates of $\beta_1,\dots,\beta_{g-l}$,
\item after surgering out the $\{\beta_1,\dots,\beta_{g-l}\}$, the
  induced curves $\beta_i$ and $\gamma_i$ (for $i=g-l+1,\dots,g$)
  lie in the punctured torus $F_i$.
\item for $i=g-l+1,\dots,g$, the curves $\beta_i$ represent meridians for the link components, disjoint from all $\gamma_j$ for $i\neq j$, and meeting $\gamma_i$ in a single transverse intersection point.
\item for $i=g-l+1,\dots,g$, the homology classes of the $\gamma_i$ correspond to the framings of the link components.
\end{itemize}

A related construction is as follows. Given $Y_1$ and $L$ as above,
choose a multi-pointed Heegaard diagram $(\Sigma_L, \alphas_L,
\betas_L, \zees_L, \wees_L)$ for $L \subset Y_1$ as in \cite{Links},
of some genus $g$; here, $\alphas_L=\{\alpha_1,\dots,\alpha_{g+l-1}\}$
and $\betas_L=\{\beta_1,\dots,\beta_{g+l-1}\}$.  Precisely, the sets
$\zees = \{z_1, \dots, z_l\}$ and $\wees = \{w_1, \dots, w_l\}$ are
collections of distinct points on $\Sigma$ disjoint from the alpha and
the beta curves, with the following two properties: first, each
connected component of $\Sigma_L\setminus\alphas_L$ and
$\Sigma_L\setminus\betas_L$ contains a single $z_i$ and a
corresponding $w_i$.  Second, if $c$ is an $l$-tuple of embedded arcs
in $\Sigma_L\setminus\betas_L$ connecting $z_i$ to $w_i$
($i=1,\dots,l$), and $c'$ is an $l$-tuple of embedded arcs in
$\Sigma_L\setminus\alphas_L$ connecting $z_i$ to $w_i$ ($i=1,\dots,l$)
then the link $L$ is the union of small push offs of $c$ and $c'$ into
the two handlebodies (induced by the beta and alpha curves,
respectively). Next, we attach handles $h_i$ to $\Sigma_L$ connecting
$z_i$ to $w_i$ (for $i=1,\dots,l$), and obtain a new surface $\Sigma$.
We choose a new $\beta_i$ ($i=g+l,\dots,g+2l-1$) to be the belt circle
of the handle $h_{i-g-l+1}$, and a new $\alpha_i$ to be the union of
the core of $h_{i-g-l+1}$ with $c'_i$, so $\alpha_i$ intersects
$\beta_i$ in one point.  Choose $\gamma_i$ ($i=1,\dots,g+l-1$) to be a
small isotopic translate of $\beta_i$, intersecting $\beta_i$ in two
points.  Let $\gamma_i^0$ ($i=g+l,\dots,g+2l-1$) be the union of $c_i$
with a core of the handle $h_{i-g-l+1}$. Obtain $\gamma_i$
($i=g+l,\dots,g+2l-1$) by applying Dehn twists to $\gamma_i^0$ around
$\beta_i$; the framing of the link is determined by the number of Dehn
twists.

In this fashion we obtain a multi-pointed triple Heegaard diagram
$(\Sigma, \alphas, \betas, \gammas, \zees),$ with $g+2l-1$ curves of
each kind. Note that $(\Sigma,\alphas,\betas,\zees)$ is an $l$-pointed
Heegaard diagram for $Y_1$, $(\Sigma,\alphas,\gammas,\zees)$ is a
$l$-pointed Heegaard diagram for $Y_2$, and
$(\Sigma,\betas,\gammas,\zees)$ is a $l$-pointed Heegaard diagram for
$\#^{g}(S^1\times S^2)$.

The new circles $\alpha_i$, $i=g+l,\dots,g+2l-1$, are part of a
maximal homologically linearly independent subset of $\alphas$, and
similarly for $\beta_i$ and $\gamma_i$ ($i=g+l,\dots,g+2l-1$).
Consequently, by the proof of Lemma~\ref{lemma:splitDiagrams}, there
is a split diagram strongly equivalent to $(\Sigma, \alphas, \betas,
\gammas, \zees),$ whose reduction
$(\Sigma',\alphas',\betas',\gammas',z')$ is obtained from $(\Sigma,
\alphas, \betas, \gammas, \zees)$ by forgetting $l-1$ of the
$\alpha_i$ (respectively $\beta_i$, $\gamma_i$), $1\leq i\leq g+l-1$,
as well as $z_2,\dots,z_k$.  It is then not hard to see that
$(\Sigma',\alphas',\betas',\gammas',z')$ is a triple Heegaard diagram
subordinate to a bouquet for $L$.

In this section we will show that for any two-handle addition, one can
construct a \emph{nice} triple Heegaard diagram strongly equivalent to
a diagram $(\Sigma, \alphas, \betas, \gammas, \zees)$ as above.  This
involves finessing the diagram for the link $L$ inside $Y$ and then,
after adding the handles and the new curves, modifying the diagram in
several steps to make it nice. For the most part we will focus on the
case when we add a single two-handle. In the last subsection we will
explain how the arguments generalize to several two-handles.

To keep language concise, in this section we will refer to elementary
domains as {\em regions}.

\subsection{A property of nice Heegaard diagrams} 

Let $(\Sigma, \alphas, \betas, \zees)$ be a multi-pointed Heegaard diagram. Recall that the diagram is called nice if all unpunctured regions are either bigons or squares. 

\begin{lemma}\label{lemma:bad_region_adjacent} 
On a nice Heegaard diagram $(\Sigma, \alphas, \betas, \zees)$, for any alpha circle $\alpha_i$ with an arbitrary orientation, there exists a punctured region $D$ which contains an edge $e$ belonging to $\alpha_i$, and such that $D$ is on the left of $\alpha_i$. The same conclusion holds for each beta circle. 
\end{lemma}

\begin{proof} 
Suppose a half-neighborhood on the left of the alpha circle $\alpha_i$ is disjoint from all the punctured regions. Then immediately to the left of $\alpha_i$ we only have good regions. There are two possibilities as indicated in Figure \ref{fig:bad_region_adjacent}. 

\begin{figure}[htbp]
 \center{\includegraphics[width=350pt]{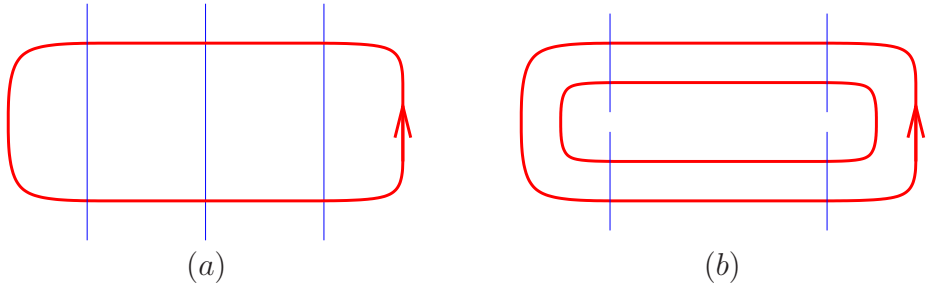}}
\caption{\label{fig:bad_region_adjacent}
{\bf Alpha curve not adjacent to the bad region.} The thick curves are alpha curves and the thin ones are beta curves.}
\end{figure}

If there is a bigon region on the left of $\alpha_i,$ then the other edge is some beta edge $\beta_j$. The region on the other side of $\beta_j$ must be a bigon region or a square since otherwise we would have a punctured region on the left of $\alpha_i$. If we reach a square, we continue to consider the next region. Eventually we will reach a bigon region since the number of regions are finite and we will not reach the same region twice. All regions involved form a disk bounded by $\alpha_i$ (as in Figure \ref{fig:bad_region_adjacent} (a)). In particular, this means $\alpha_i$ is null homologous. This contradicts the fact that the $g+k-1$ alpha circles represent linearly independent classes in $H_1(\Sigma\setminus \zees).$

In the second case, there are no bigon regions. Then on the left of $\alpha_i$, we see a chain of squares, as in Figure \ref{fig:bad_region_adjacent} (b). The opposite edges on these squares give another alpha circle, say $\alpha_j$. Then $\alpha_i$ and $\alpha_j$ are homologous to each other in $H_1(\Sigma\setminus \zees).$ This contradicts the same fact as in the previous case. 
\end{proof}

Recall that in order to define the triangle maps it is necessary for the triple Heegaard diagram to be weakly admissible in the sense of~\cite[Definition 8.8]{HolDisk}. 

\begin{corollary}\label{lemma:admissible}
If $(\Sigma,\alphas,\betas,\gammas,\zees)$ is a nice multi-pointed triple Heegaard diagram then $(\Sigma,\alphas,\betas,\gammas,\zees)$ is weakly admissible.
\end{corollary}

\begin{proof}
By definition, the diagram is weakly admissible if there are no nontrivial domains $D$ supported in $\Sigma \setminus \zees$ with nonnegative multiplicity in all regions, and whose boundary is a linear combination of alpha, beta, and gamma curves. Suppose such a domain $D$ exists, and consider a curve appearing with a nonzero multiplicity in $\partial D.$ Without loss of generality, we can assume this is an alpha curve, and all regions immediately to its left have positive multiplicity in $D.$ By Lemma~\ref{lemma:forget}, the diagram $(\Sigma,\alphas,\betas,\zees)$ is nice. Lemma~\ref{lemma:bad_region_adjacent} now gives a contradiction.
\end{proof}

\subsection{A single two-handle addition} 
\label{sec:boss} 

Let $(Y,K)$ be a three-manifold together with a knot $K \subset Y$. We choose a singly pointed Heegaard diagram $(\Sigma, \alphas, \betas, z)$ for $Y$ together with an additional basepoint $w \neq z \in \Sigma \setminus (\alphas \cup \betas)$ such that the two basepoints determine the knot as in \cite{Knots}. After applying the algorithm from \cite{SarkarWang} to the Heegaard diagram, we can assume that the Heegaard diagram is nice, with $D_z$ the (usually bad) region containing the basepoint $z.$ Furthermore, the algorithm in \cite{SarkarWang} also ensures that $D_z$ is a polygon.

We denote by $D_w$ the region containing $w;$ note that either $D_w = D_z$ or $D_w$
is good. Throughout this section, we will suppose that $D_w$ and $D_z$ are two
different regions, and that $D_w$ is a rectangle. The case when $D_w = D_z$
corresponds to surgery on the unknot, which is already well understood. The case when
$D_w$ is a bigon can be avoided by modifying the original
  diagram by a finger move. (Alternately, this case can be treated similarly to the case that $D_w$ is
  a rectangle.)

Let $W$ be the four manifold with boundary obtained from $Y\times[0,1]$ by adding a two handle along $K$ in $Y\times\{1\}$, with some framing. $W$ gives a cobordism between $Y$ and $Y'$, where $Y'$ is obtained from $Y$ by doing the corresponding surgery along $K$.

Now we are ready to describe our algorithm to get a nice triple Heegaard diagram for the cobordism $W.$

\subsection*{Step 1. Making the knot embedded in the Heegaard diagram}\

Let $c$ be an embedded arc in $\Sigma$ connecting $z$ and $w$ in the complement of beta curves, and $c'$ be an embedded arc connecting $z$ and $w$ in the complement of alpha curves. The union of $c$ and $c'$ is a projection of the knot $K \subset Y$ to the surface $\Sigma,$ where $\Sigma$ is viewed as a Heegaard surface in $Y.$ For convenience, we will always assume that $c$ and $c'$ do not pass through any bigon regions, and never leave a rectangle by the same edge through which they entered; this can easily be achieved.

In this step, we modify the doubly pointed Heegaard diagram $(\Sigma, \alphas, \betas, z,w)$ to make $c \cup c'$ embedded in $\Sigma,$ while preserving the niceness of the Heegaard diagram.

Typically, $c$ and $c'$ have many intersections. We modify the diagram inductively by stabilization at the first intersection $p\in D_p$ on $c'$ (going from $z$ to $w$) to remove that intersection, while making sure that the new diagram is still nice.

\begin{figure}[htbp]
 \center{\includegraphics[width=400pt]{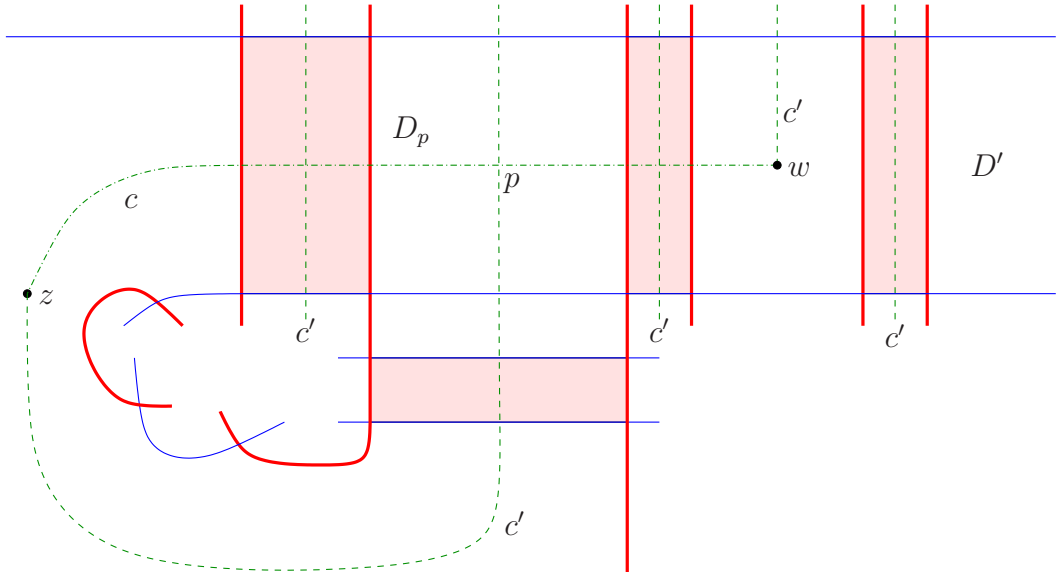}} 
\caption{\label{fig:knot_embedded_before} 
{\bf Projection of the knot to the Heegaard surface.} As usual, the thick lines are alpha curves, and the thin lines are betas. The point $p$ is the first on the dashed $c'$ curve (starting from $z$) where an intersection with $c$ takes place. The light shading in the upper three domains indicates that there might be more parallel copies of alpha curves or parts of the $c'$ curve. Similarly, the light shading in the lower domain indicates the presence of an arbitrary number of parallel beta segments there. Note that some of the regions in the rightmost shaded domain on the top can coincide with some of the regions in the bottom shaded domain; however, this fact will not create any difficulties.} 
\end{figure}

A neighborhood of $c$ and the part on $c'$ from $z$ to $p$ are shown in Figure \ref{fig:knot_embedded_before}. In the same picture, if we continue the chain of rectangles containing $c$, we will end up with a region $D^\prime$ which is either a bigon or the punctured region $D_z$. 

\begin{figure}[htbp]
 \center{\includegraphics[width=400pt]{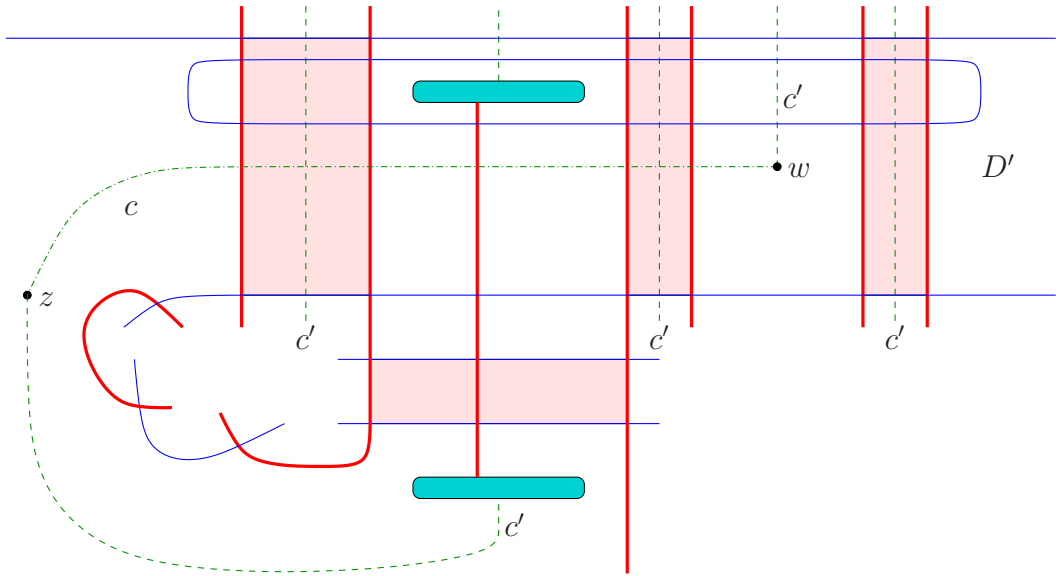}} 
\caption{\label{fig:knot_embedded_after} 
{\bf Stabilizing at the intersection.} The two darkly shaded ovals are the two feet of the handled we added. In the stabilized picture, we stretch the new beta curve until it reaches the regions $D_z$ and $D'$, so that the diagram is still nice.} 
\end {figure}

To get rid of the intersection point $p$, we stabilize the diagram as in Figure \ref{fig:knot_embedded_after}. More precisely, we do a stabilization followed some handleslides of the beta curves and an isotopy of the new beta curve. After these moves, the number of intersection points decreases by one and the diagram is still nice. 

If we iterate this process, in the end we get a nice Heegaard diagram in which $c$ and $c'$ only intersect at their endpoints. Furthermore, the bad region $D_z$ is still a polygon.

\subsection*{Step 2. Adding twin gamma curves.} \

Our goal in Steps 2 and 3 is to describe a particular triple Heegaard diagram for the cobordism $W$. Starting with the alpha and the beta curves we already have, for each beta curve $\beta_i$ we will add a gamma curve $\gamma_i$ (called its twin) which is isotopic to $\beta_i$ and intersects it in exactly two points. (After this, we will add some more curves in the next step.)
 
\begin{figure}[htbp]
 \center{\includegraphics[width=350pt]{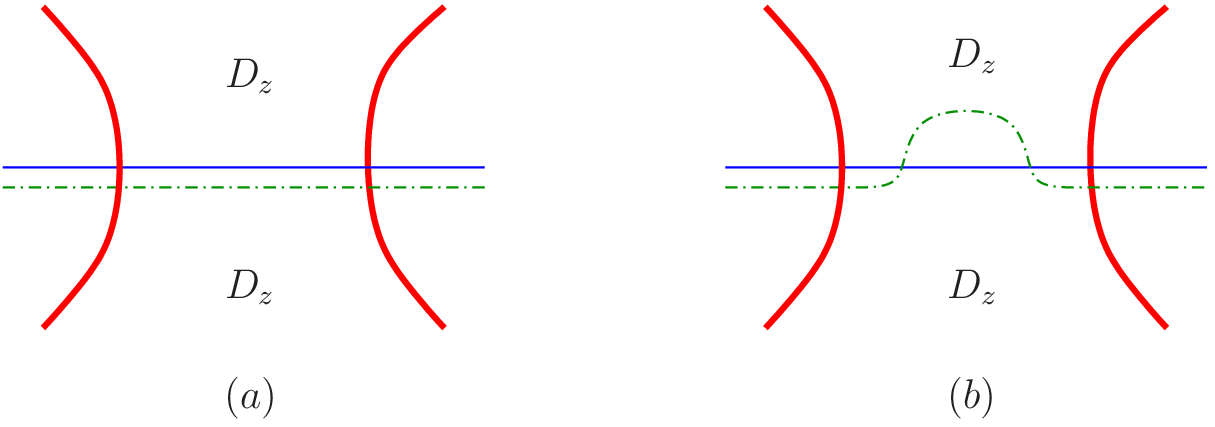}} 
\caption{\label{fig:adding_gamma_curve_simple}
{\bf Adding twin gamma curves - Case 1.}} {Here and after, without further specification, we make the convention that the thick arcs are alpha arcs, the thin ones are beta arcs, and the interrupted ones are gamma arcs.}
\end{figure}

For any beta curve $\beta_i$, by Lemma \ref{lemma:bad_region_adjacent} we can choose a region $D_i$ so that $D_i$ is adjacent to the punctured region $D_z$ with the common edge on $\beta_i$. If $D_i=D_z$, then we add 
$\gamma_i$ close and parallel to $\beta_i$ as in Figure \ref{fig:adding_gamma_curve_simple} (a), and make a finger move as in Figure \ref{fig:adding_gamma_curve_simple} (b).

Suppose now that $D_i$ is different from $D_z$. Then $D_i$ has to be a good region.

\begin{figure}[htbp]
 \center{\includegraphics[width=350pt]{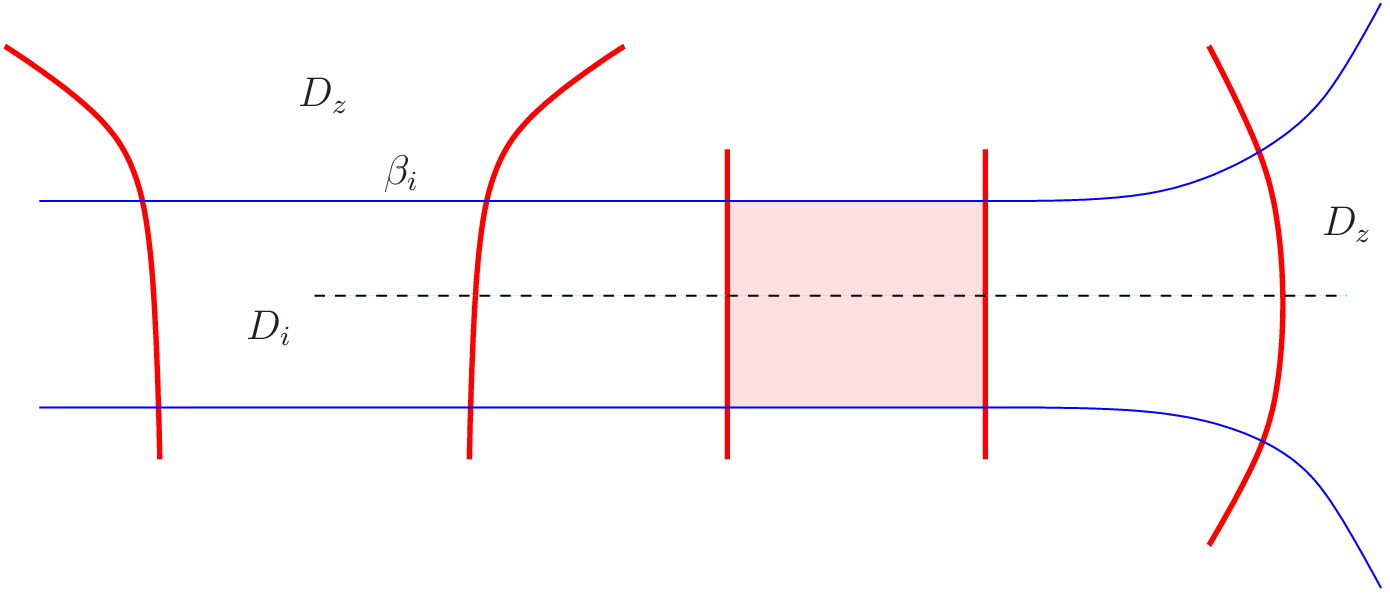}} 
\caption{\label{fig:adding_gamma_curve_before}
{\bf An arc connecting $D_i$ to $D_z$.} The shaded domain may contain several parallel alpha edges, and the dashed arc is the connecting arc.}
\end{figure}

\begin{figure}[htbp]
 \center{\includegraphics[width=350pt]{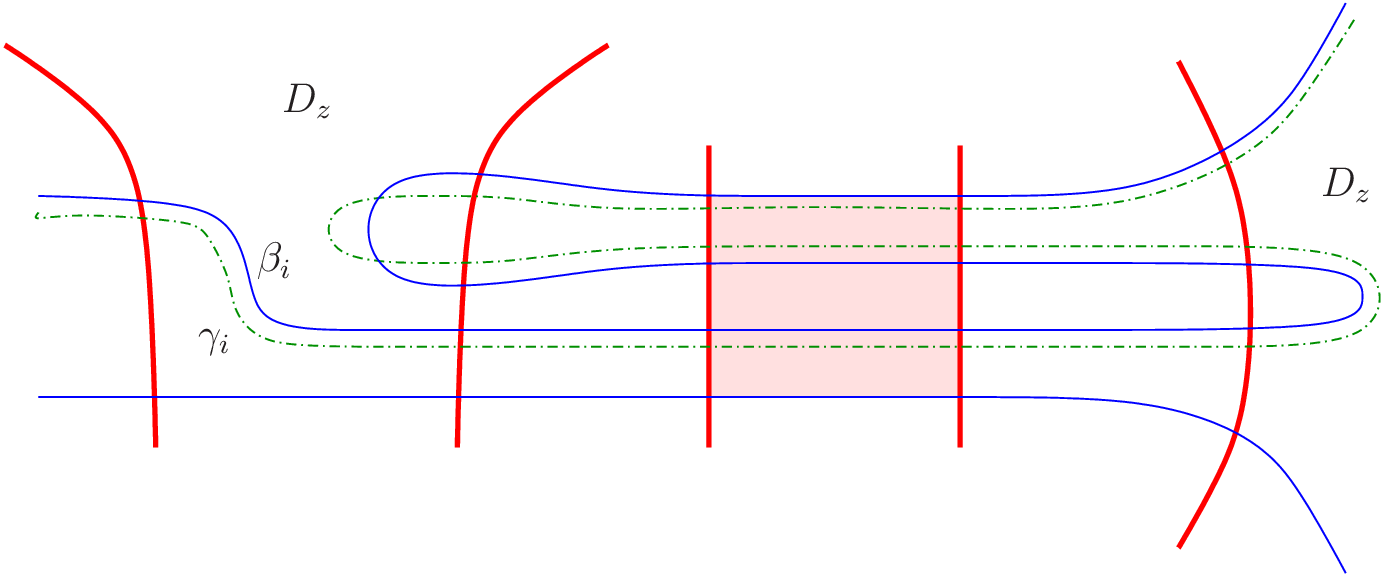}} 
\caption{\label{fig:adding_gamma_curve_after}
{\bf Adding twin gamma curves - Case 2.} Before adding $\gamma_i,$ we do the finger move shown here. Again, the shaded domain may contain several parallel alpha segments.}
\end{figure}

If $D_i$ is not a bigon, since the complement of the beta curves in $\Sigma$ is connected, we can connect $D_i$ with $D_z$ without intersecting beta curves, via an arc traversing a chain of rectangles, as indicated in the Figure 
\ref{fig:adding_gamma_curve_before}. Then we do a finger move of the curve $\beta_i$ as indicated in Figure \ref{fig:adding_gamma_curve_after}. Note that the knot remains embedded in $\Sigma$.

Now we have a bigon region. We then add the gamma curve $\gamma_i$ as shown in Figure \ref{fig:adding_gamma_curve_after}.

\begin{figure}[htbp]
 \center{\includegraphics[width=200pt]{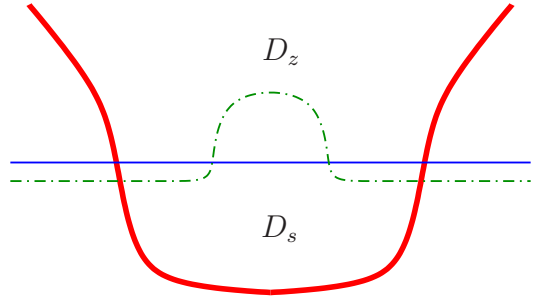}}
\caption{\label{fig:beta_gamma_bigon}
{\bf Bigon between beta and gamma curves.} The region $D_s$ will be dealt with in a special way in later steps.}
\end{figure}

Note that for each pair $\beta_i$ and $\gamma_i$, we either have one sub-diagram of the form in Figure~\ref{fig:adding_gamma_curve_simple} (b), or one sub-diagram of the form in Figure \ref{fig:beta_gamma_bigon}. Observe also that during this process, no bad region other than $D_z$ is created.

\subsection*{Step 3. Stabilization and two-handle addition}\

\begin{figure}[htbp]
 \center{\includegraphics[width=400pt]{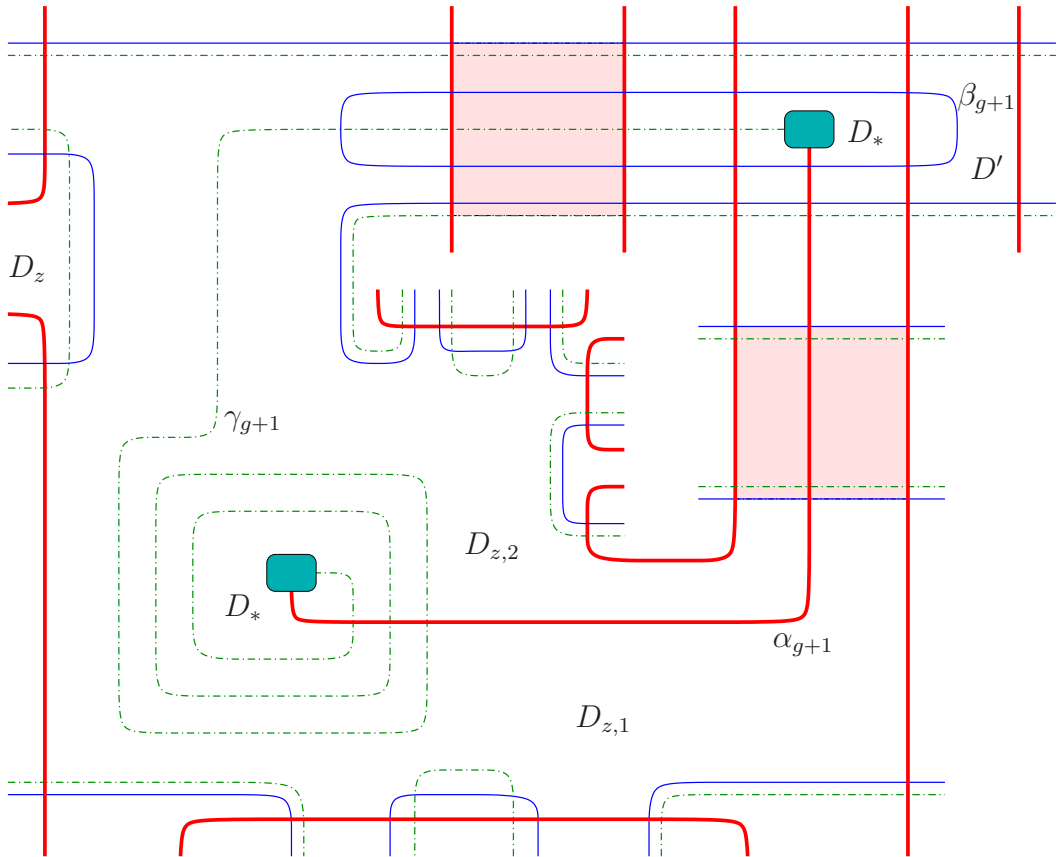}}
\caption{\label{fig:stabilization}
{\bf Stabilization - Case 1.} The two small darkly shaded ovals are the two feet of the handle. The upper large lightly shaded area may have several parallel copies of alpha arcs, while the lower one may have several copies of beta and gamma arcs. The upper left $D_z$ denotes  either $D_{z,1}$ or $D_{z,2};$ such a domain can occur at various places on the boundaries of $D_{z,1}$ or $D_{z,2},$ and it corresponds to Case 1 in Step 2.} 
\end{figure}

After Step 2, the knot is still embedded in the Heegaard diagram. In other words, we can use arcs to connect $z$ to $w$ by paths in the complement of alpha curves, and in the complement of beta curves so that the two arcs do not intersect except for the end points $w$ and $z$, and do not pass through any bigons. We will see two chains of squares, as indicated in Figure \ref{fig:stabilization}.

We do a stabilization of the Heegaard diagram by adding a handle with one foot in each of $D_z$ and $D_w$. We add the additional beta circle $\beta_{g+1}$ to be the meridian of the handle, which we push along $c$ until it reaches $D_z.$ We also push $\beta_{g+1}$ through the opposite alpha edge of $D_w$, into the adjacent region. Then, we connect the two feet in the complement of alpha curves along $c'$ and get a new alpha circle $\alpha_{g+1}$. Finally, we add the surgery gamma circle $\gamma_{g+1}$ as in Figure \ref{fig:stabilization}. The result is a triple Heegaard diagram (with $3(g+1)$ curves) which represents surgery along the knot $K \subset Y$, with a particular framing; the framing is the sum of the number of twists of $\gamma_{g+1}$ around the handle and a constant depending only on the original Heegaard diagram.

\begin{figure}[htbp]
 \center{\includegraphics[width=350pt]{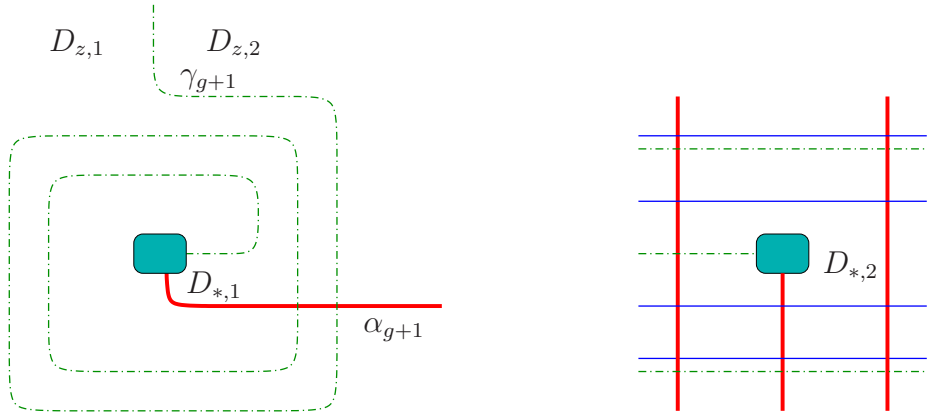}}
\caption{\label{fig:stabilization_anti}
{\bf Stabilization - Case 2.} If the $\gamma_{g+1}$ twists around the handle in the opposite direction, we still have a picture similar to Figure~\ref{fig:stabilization}. The only differences are in the neighborhoods of the two feet of the handle, which are shown here. }
\end{figure}

Note that, depending on the framing, the local picture around the two feet of the handle may also look like Figure \ref{fig:stabilization_anti}, in which case instead of the octagon region $D_*$ from Figure \ref{fig:stabilization} we have two hexagon bad regions $D_{*,1}$ and $D_{*,2}$.

After the stabilization, $\alpha_{g+1}$ and $\gamma_{g+1}$ separate $D_z$ into several regions; among these, $D_{z,1}$ and $D_{z,2}$ are (possibly) bad but all other regions are good. We end up with a diagram with four (or five) bad regions: $D_{z,1}$, $D_{z,2}$, $D_*$ (or $D_{*,1}$ and $D_{*,2}$), and $D^\prime$. (In some cases, $D_{z,1}$ or $D_{z,2}$ might be good, or, if there is little winding of $\gamma_{g+1},$ some of $D_*$, $D_{z,1},$ and $D_{z,2}$ might coincide. The argument in these cases is a simple adaptation of the one we give below.) We will kill the badness of $D^\prime,$ $D_{z,2},$ and $D_*$ (or $D_{*,1}$ and $D_{*,2}$), while the region $D_{z,1}$ will be the one containing the basepoint $z$ for our final triple Heegaard diagram.

\subsection*{Step 4. Killing the bad region $D^\prime$}\ 

We push the finger in $D^\prime$ across the opposite alpha edge until we reach a bigon, $D_{z,1}, D_{z,2},$ or a region of type $D_{s}$ as in Figure \ref{fig:beta_gamma_bigon}.

\subsubsection*{Case 1. A bigon is reached}

\begin{figure}[htbp]
 \center{\includegraphics[width=400pt]{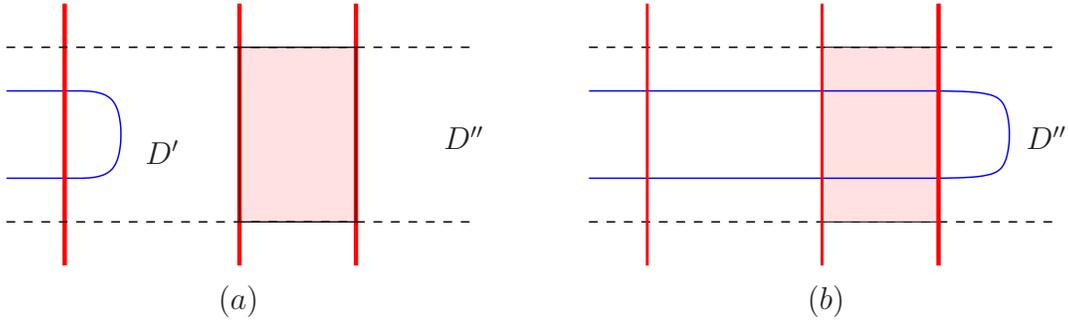}}
\caption{\label{fig:dprime_bigon}
{{\bf Killing the badness of $D^\prime$: bigon, $D_{z,1}$ or $D_{z,2}$ case.}} The thick arcs are alpha curves, the thin arcs represent $\beta_{g+1}$, and the dashed arcs can be either betas or gammas. The shaded part has several parallel alpha arcs. The rightmost region $D''$ is a bigon, $D_{z,1},$ or $D_{z,2}.$} 
\end{figure}

In this case (Figure \ref{fig:dprime_bigon} (a)), our finger move will kill the badness of $D^\prime$, as indicated in Figure \ref{fig:dprime_bigon} (b), and does not create any new bad regions.

\subsubsection*{Case 2. $D_{z,1}$ or $D_{z,2}$ is reached}
This is completely similar to Case 1. The finger move kills the badness of $D^\prime$, and does not create any new bad regions.

\subsubsection*{Case 3. A region of type $D_s$ is reached}

Let us suppose the topmost region in Figure \ref{fig:beta_gamma_bigon}
is $D_{z,1}$. The case when the topmost region is $D_{z,2}$ is
completely similar.

\begin{figure}[htbp]
 \center{\includegraphics[width=350pt]{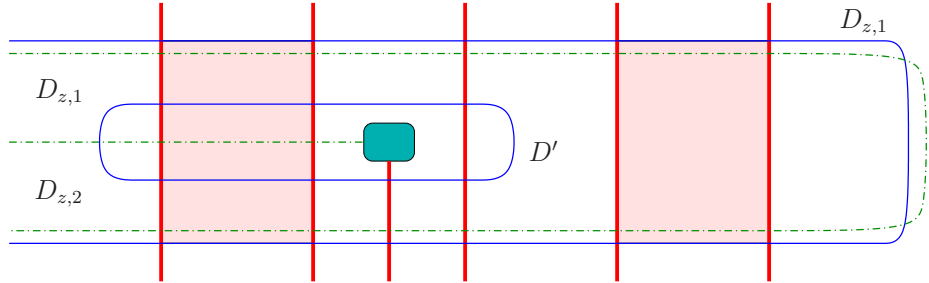}}
\caption{\label{fig:dprime_ds_bad}
{{\bf Killing the badness of $D^\prime$: $D_s$ type region - before.}} The smaller shaded region is the foot of the handle inside $D_w$, while the larger two shaded regions may contain several parallel alpha arcs.} 
\end{figure}

\begin{figure}[htbp]
 \center{\includegraphics[width=350pt]{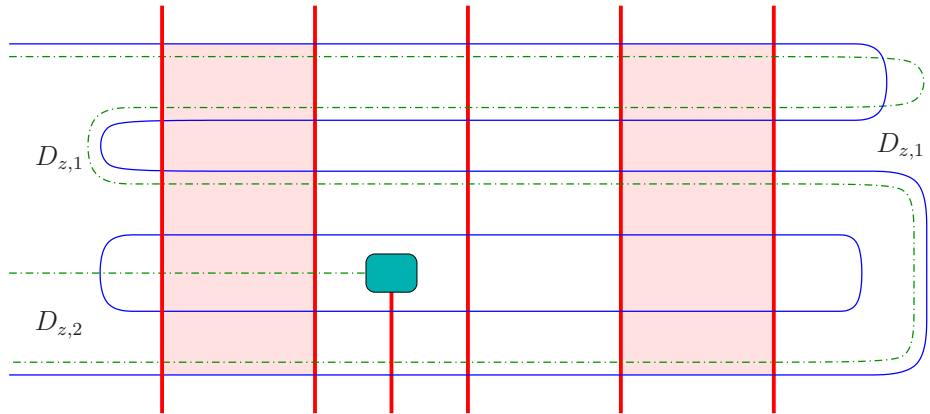}} 
\caption{\label{fig:dprime_ds_nice} 
{\bf Killing the badness of $D^\prime$: $D_s$ type region - after.} Conventions and shaded regions are the same as in Figure \ref{fig:dprime_ds_bad}. This figure is obtained from 
Figure \ref{fig:dprime_ds_bad} via a double finger move.} 
\end{figure}

The regions involved look like Figure \ref{fig:dprime_ds_bad}. If on the left $D_{z,1}$ is on top of $D_{z,2}$, we isotope the diagram to look as in Figure \ref{fig:dprime_ds_nice}. The case when on the left of Figure \ref{fig:beta_gamma_bigon} $D_{z,2}$ is on top of $D_{z,1}$ is similar, except that we do the double finger move on the other side of $\beta_{g+1}$.

We have now killed the badness of $D^\prime$.

\subsection*{Step 5. Killing the badness of $D_{z,2}$}\ 

\begin{figure}[htbp]
 \center{\includegraphics[width=400pt]{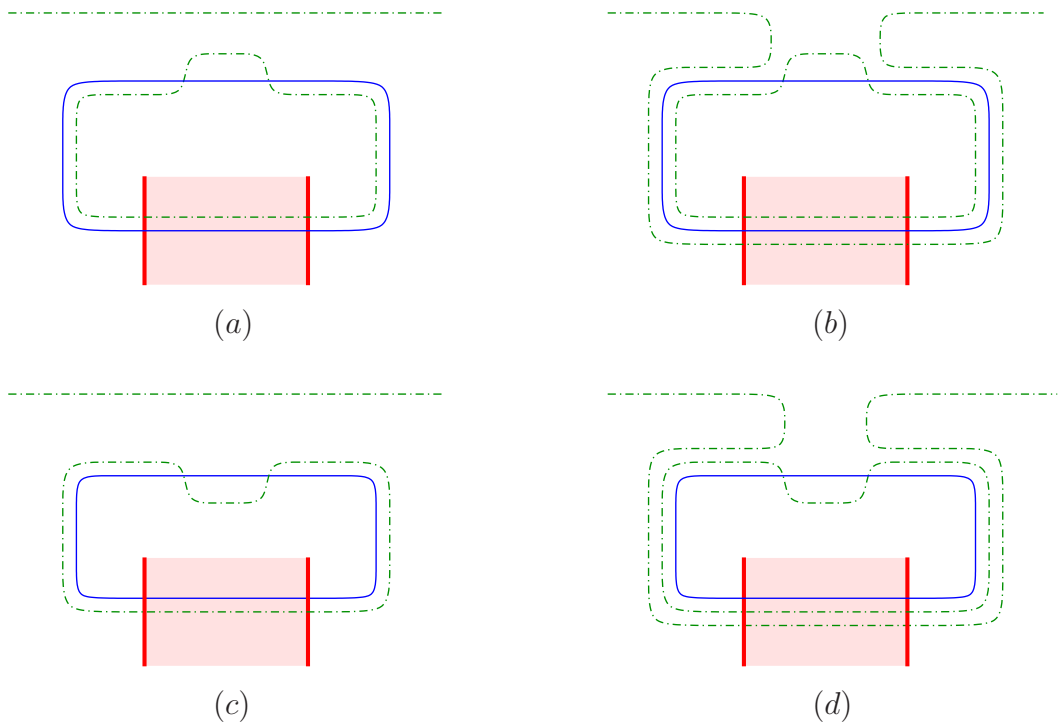}}
\caption{\label{fig:dw2_special_handleslide}
{\bf Killing the badness of $D_{z,2}$: special handleslides.} There might be more alpha arcs in the shaded regions.} 
\end{figure}

If there are any bigons between beta and gamma curves adjacent to $D_{z,2}$ as in Figure \ref{fig:beta_gamma_bigon} or Figure \ref{fig:adding_gamma_curve_simple} (b), also shown in Figures \ref{fig:dw2_special_handleslide} (a) resp. (c), we do a ``handleslide'' (more precisely, a handleslide followed by an isotopy) of $\gamma_{g+1}$ over each $\gamma_i$ ($i\leq g$) involved as indicated in Figures \ref{fig:dw2_special_handleslide} (b) resp. (d).

\begin{figure}[htbp]
 \center{\includegraphics[width=440pt]{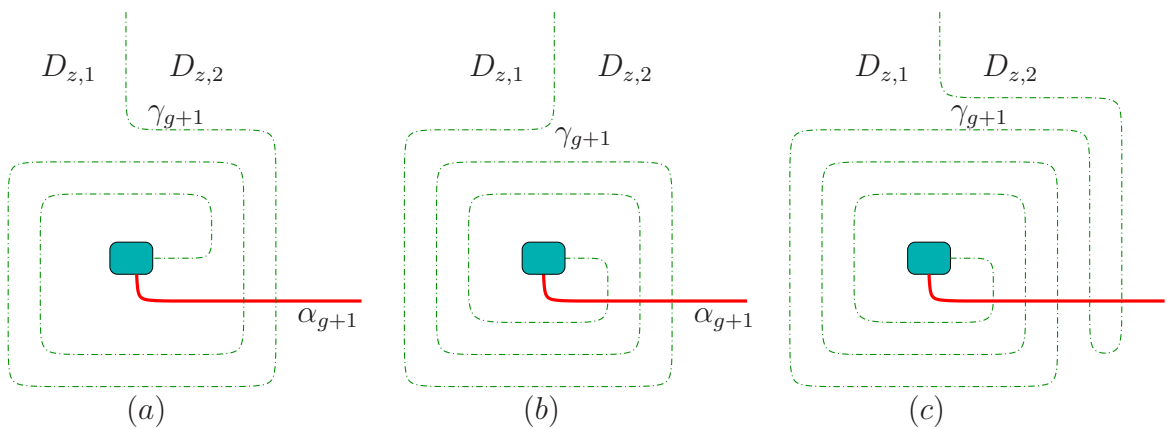}}
\caption{\label{fig:dw2_two_patterns}
{\bf Killing the badness of $D_{z,2}$: two patterns.} The curve $\gamma_{g+1}$ can rotate around the shaded oval either as in (a) or as in (b). If (b) occurs, we replace it with (c).} 
\end{figure}

The intersection of $\gamma_{g+1}$ and $\alpha_{g+1}$ has the pattern as Figure \ref{fig:dw2_two_patterns} (a) or (b). In case (a), we do nothing. In case (b), we do the finger move as in Figure \ref{fig:dw2_two_patterns} (c).

\begin{figure}[htbp]
 \center{\includegraphics[width=350pt]{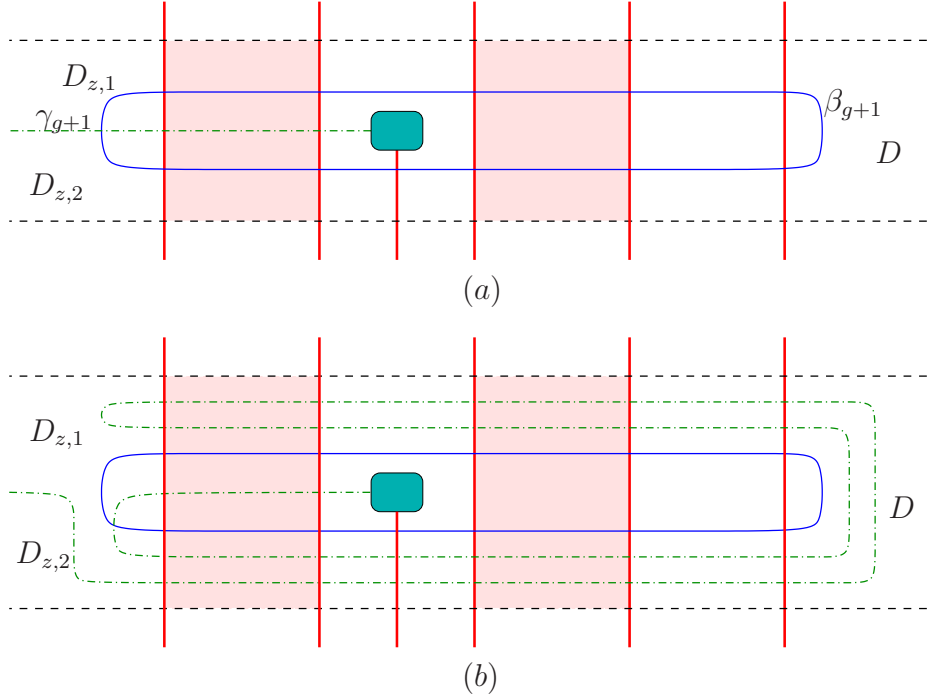}}
\caption{\label{fig:dw2_beta_gamma_crossing}
{\bf Killing the badness of $D_{z,2}$: special beta-gamma crossing.} The shaded area may contain several parallel alpha arcs. The dashed arcs can be beta arcs or gamma arcs depending on different cases. The region $D$ is a bigon, $D_{z,1},$ or $D_{z,2}$. We isotope (a) into (b) in order to remove the beta-gamma intersection point on the boundary of $D_{z,2}.$} 
\end{figure}

Now among the possibly bad regions generated from $D_{z,2}$, we have a unique one whose boundary has an intersection of a beta curve with a gamma curve, namely the one near $\beta_{g+1}$ as in Figure \ref{fig:dw2_beta_gamma_crossing} (a). (See also Figures \ref{fig:stabilization} and \ref{fig:dprime_ds_nice}.)

We then do a finger move as in Figure \ref{fig:dw2_beta_gamma_crossing} (b). Note that this finger move will not create any badness other than that of $D_{z,1}$.

After these special handleslides and finger moves, the region $D_{z,2}$ is divided into several possibly bad regions $R_1, \dots, R_m.$ These bad regions are all adjacent to $D_{z,1}$ via arcs on $\gamma_{g+1}$ and, 
furthermore, there are no intersection points of beta and gamma curves on their boundaries. We seek to kill the badness of $R_1, \dots, R_m$ using the algorithm in \cite{SarkarWang}. The algorithm there consisted of inductively decreasing a complexity function defined using the unpunctured bad regions. In our situation, we apply a simple modification of the algorithm to the Heegaard diagram made of the alpha and the gamma curves; the modification consists of the fact that we do not deal with the bad region(s) $D_*$ (or $D_{*,1}$ and $D_{*,2}$), but rather only seek to eliminate the badness of $R_1, \dots, R_m;$ thus, in the complexity function we do not include terms that involve the badness and distance of $D_*$ (or $D_{*,1}$ and $D_{*,2}$).

Since all the $R_i$'s are adjacent to the preferred (punctured) region $D_{z,1}$ via arcs on $\gamma_{g+1},$ the algorithm in \cite{SarkarWang} prescribes doing finger moves of $\gamma_{g+1}$ through alpha curves, and 
(possibly) handleslides of $\gamma_{g+1}$ over other gamma curves. We do all these moves in such a way as not to tamper with the arrangements of twin beta-gamma curves, i.e. as not to introduce any new intersection points between $\gamma_{g+1}$ and $\beta_i,$ for any $i \leq g$. (In other words, we can think of fattening $\gamma_1, \dots, \gamma_g$ before applying the algorithm, so that they include their respective twin beta curves.) In particular, regions of type $D_{s}$ are treated as bigons. 

The fact that the algorithm in \cite{SarkarWang} can be applied in this fashion is based on the following two observations: 
\begin{itemize}
 \item Our fingers or handleslides will not pass through the regions adjacent to $\beta_{g+1}$, except possibly $D_{z,1}$ itself. (This is one benefit of the modification performed in Figure~\ref{fig:dw2_beta_gamma_crossing}.)
 \item We will not reach any squares between $\beta_i$ and $\gamma_i \ (i \leq g),$ nor the ``narrow'' squares created in 
Figure \ref{fig:dw2_special_handleslide}.
\end{itemize}

\begin{figure}[htbp]
 \center{\includegraphics[width=400pt]{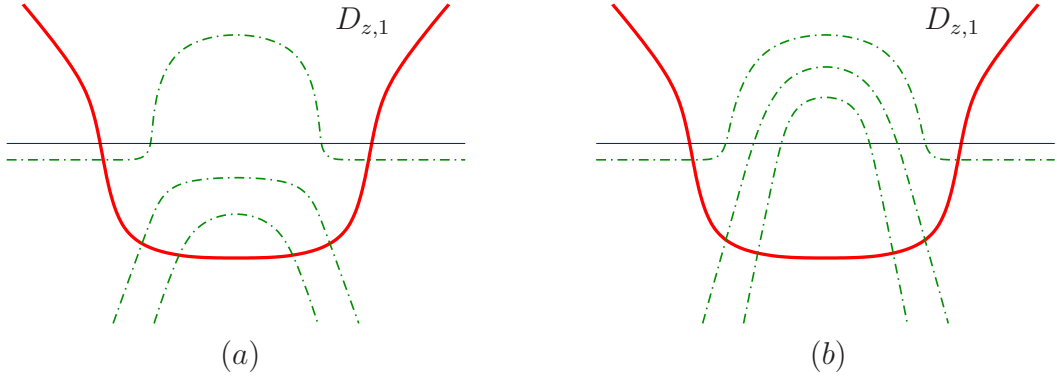}}
\caption{\label{fig:dw2_special_bad_region}
{\bf Killing the badness of special bad regions.} We isotope the gamma curves to kill the hexagon in (a).}
\end{figure}

In the end, all the badness of $R_1, \dots, R_m$ is killed. We arrive at a Heegaard diagram which might still have some bad regions coming from regions of type $D_s,$ as in Figure \ref{fig:dw2_special_bad_region} (a). We kill these bad regions using the finger moves indicated in Figure \ref{fig:dw2_special_bad_region} (b).

After these moves, the only remaining bad regions are $D_*$ (or $D_{*,1}$ and $D_{*,2}$), and the preferred bad region $D_{z,1}.$

\subsection*{Step 6. Killing the badness of $D_*$ (or $D_{*,i}$)} \

Our remaining task is to kill the badness of $D_*$ or $D_{*,i}$. Recall that depending on the pattern of the intersection of $\alpha_{g+1}$ and $\gamma_{g+1}$ (cf. Figure \ref{fig:dw2_two_patterns}), there are two cases: either we have an octagon bad region $D_*$, or two hexagon bad regions $D_{*,1}$ and $D_{*,2}$.

\begin{figure}[htbp]
 \center{\includegraphics[width=440pt]{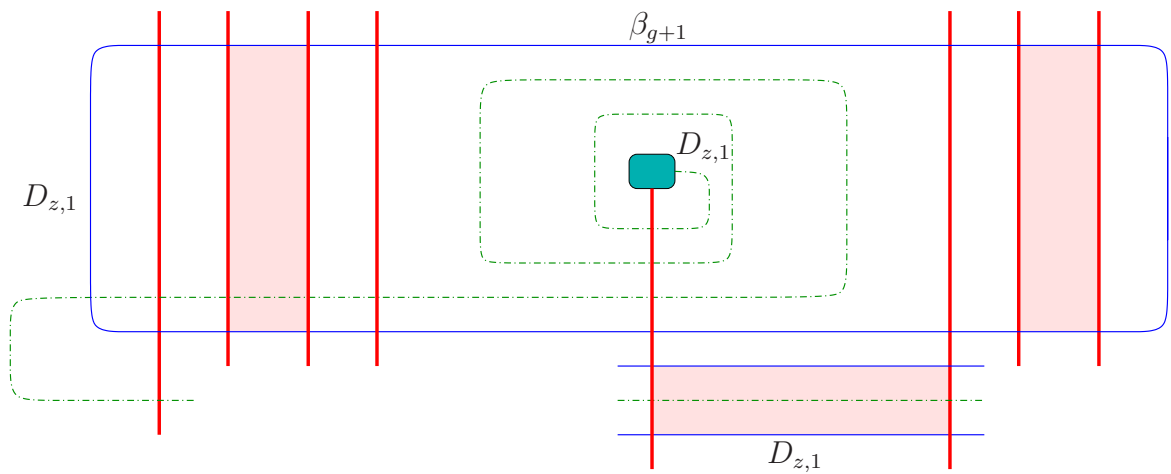}}
\caption{\label{fig:dstar_bad}
{\bf The bad region $D_*$.} Conventions are as before. Lightly shaded regions mean several parallel arcs.}
\end{figure}

\begin{figure}[htbp]
\center{\includegraphics[width=440pt]{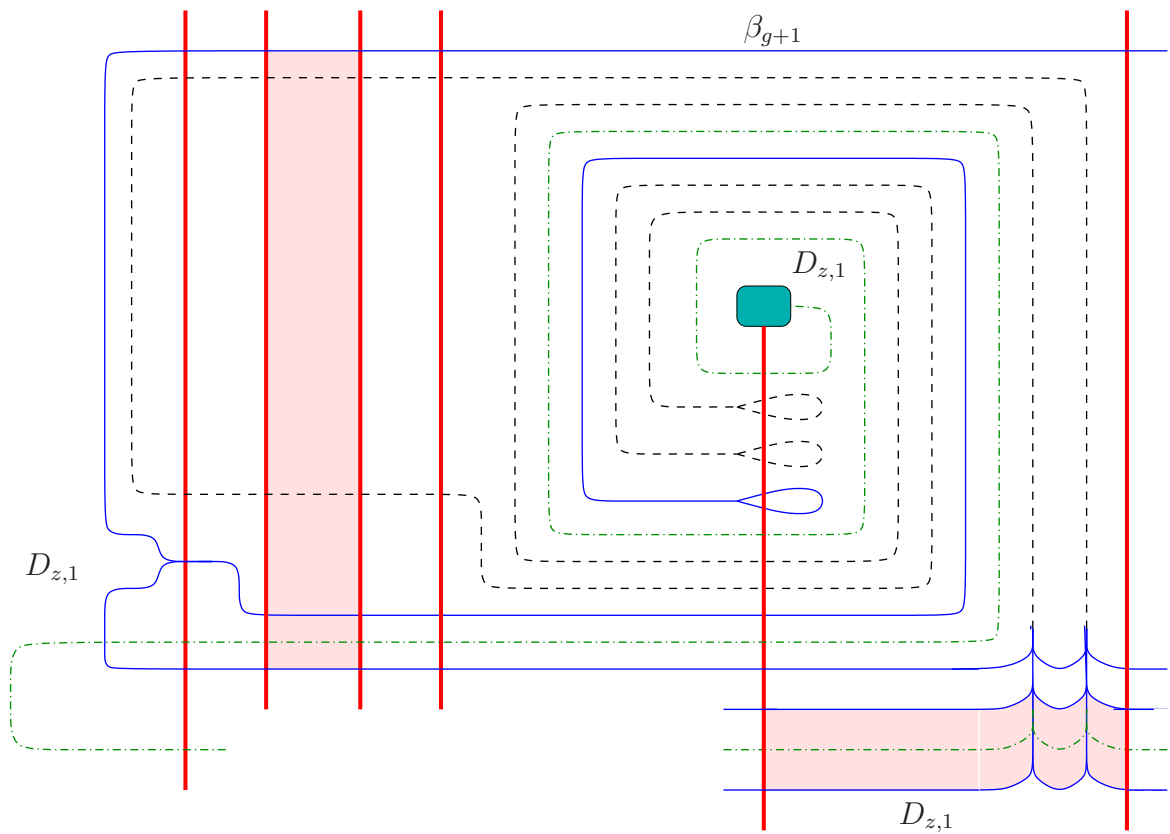}}
\caption{\label{fig:dstar_nice}
{\bf Killing the badness of $D_*$.} We push two multiple fingers (containing several beta and gamma curves) from the bottom right of the diagram, and a single finger from the left. We are using the train-track convention, cf. Figure~\ref{fig:train-track} below.} 
\end{figure}

In the first case, one possibility is that a neighborhood of $\alpha_{g+1}\cup\beta_{g+1}\cup\gamma_{g+1}$ looks as in Figure \ref{fig:dstar_bad}. We then do the finger moves indicated in Figure \ref{fig:dstar_nice}. It is routine to check that the new diagram is isotopic to the one in Figure \ref{fig:dstar_bad}. 


\begin{figure}[htbp]
 \center{\includegraphics[width=350pt]{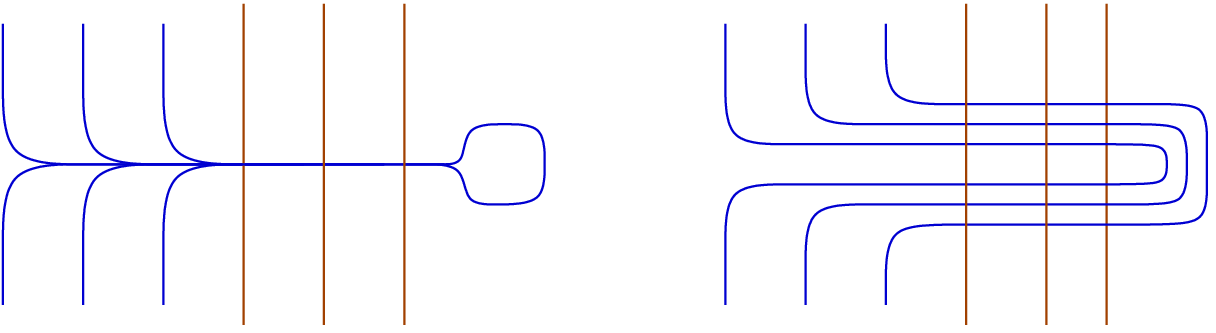}}
\caption{\label{fig:train-track} 
{\bf The train-track convention.} We use the left diagram to denote a multiple finger, i.e. the situation pictured on the right. The curves involved can be of various kinds.}
\end{figure}

\begin{figure}[htbp]
 \center{\includegraphics[width=440pt]{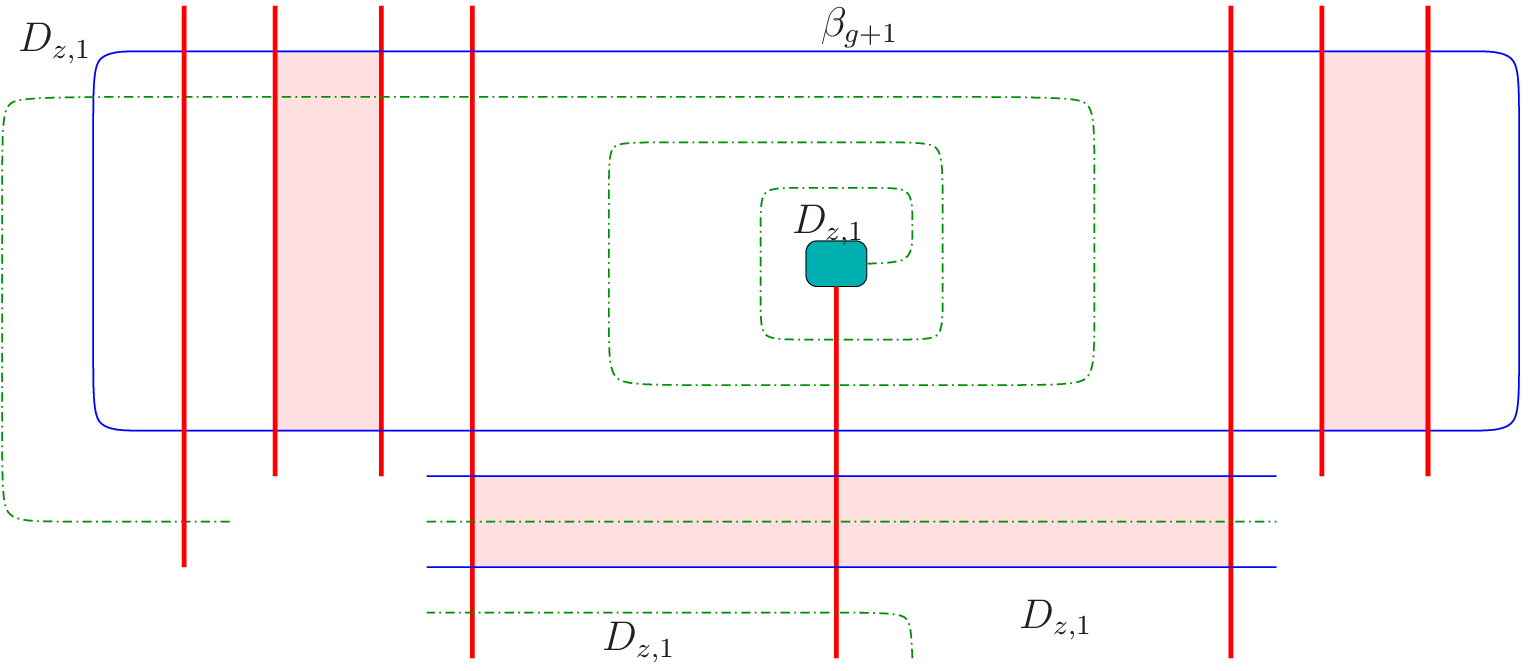}}
\caption{\label{fig:dstar_bad_case2}
{\bf The bad regions $D_{*,1}$ and $D_{*,2}$.} Going from each of the two hexagons down through beta and gamma curves, we eventually reach $D_{z,1}.$ The path from one of the two hexagons (in the picture, the one on the left) encounters one additional beta curve before reaching $D_{z,1}.$ }
\end{figure}

\begin{figure}[htbp]
 \center{\includegraphics[width=440pt]{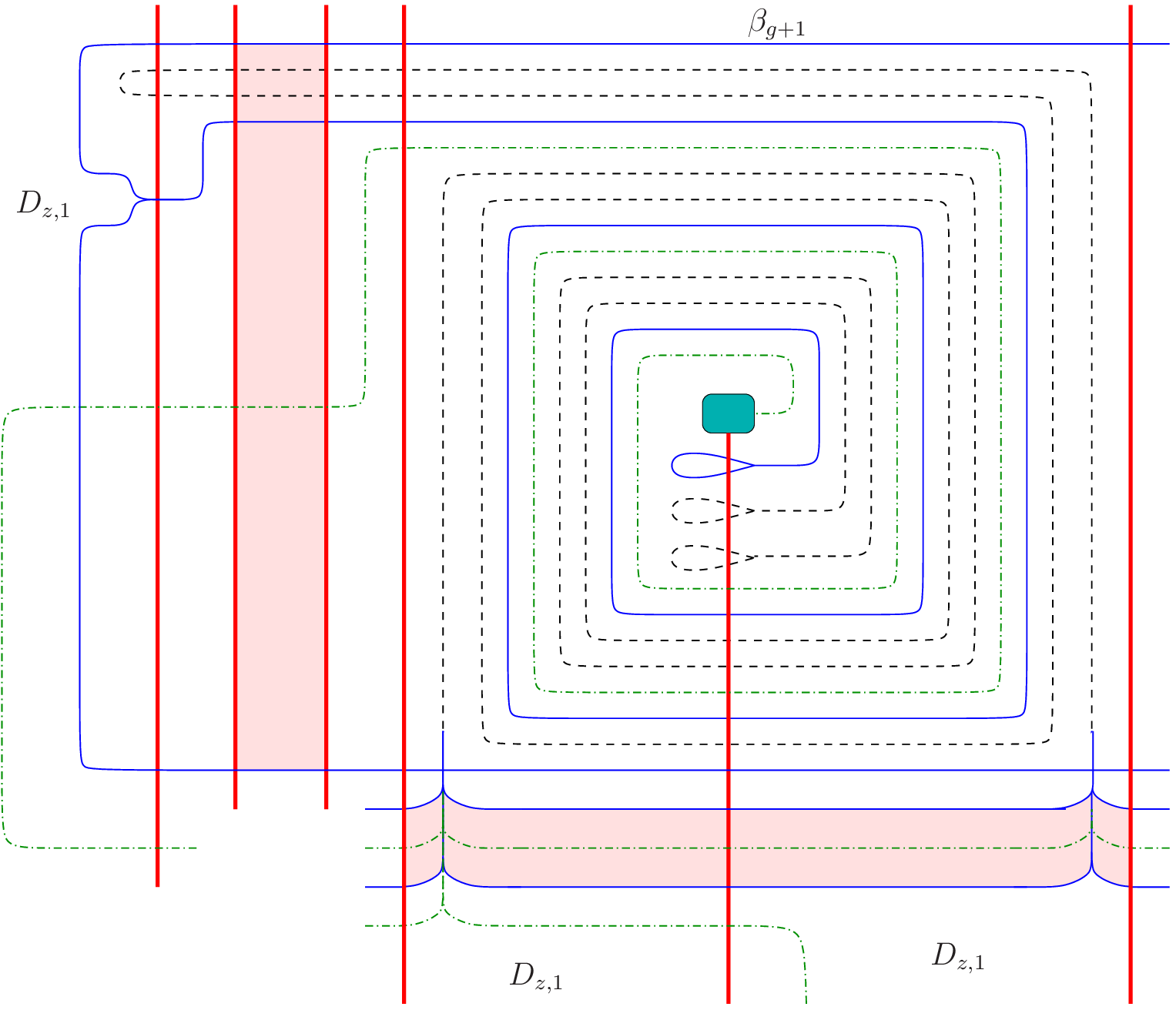}}
\caption{\label{fig:dstar_nice_case2}
{\bf Killing the badness of $D_{*,1}$ and $D_{*,2}$.} In this case, we again push two multiple fingers from the bottom, but starting from two different sides of $\alpha_{g+1}$. As before, we also push a finger form the left.} 
\end{figure}

Similarly, in the second case, one possibility is that a neighborhood of $\alpha_{g+1}\cup\beta_{g+1}\cup\gamma_{g+1}$ looks as in Figure \ref{fig:dstar_bad_case2}. In this case, we do the finger moves indicated in Figure \ref{fig:dstar_nice_case2}.

However, the actual picture on the Heegaard diagram may differ from Figure~\ref{fig:dstar_bad} or \ref{fig:dstar_bad_case2} in several (non-essential) ways.

One possible difference is that at the bottom of the Figure~\ref{fig:dstar_bad_case2}, the extra gamma curve on top of $D_{z,1}$ might be on the right rather than on the left; however, we can still push the two fingers starting from $D_{z,1}$ on each side of $\alpha_{g+1}$. 

\begin{figure}[htbp]
 \center{\includegraphics[width=440pt]{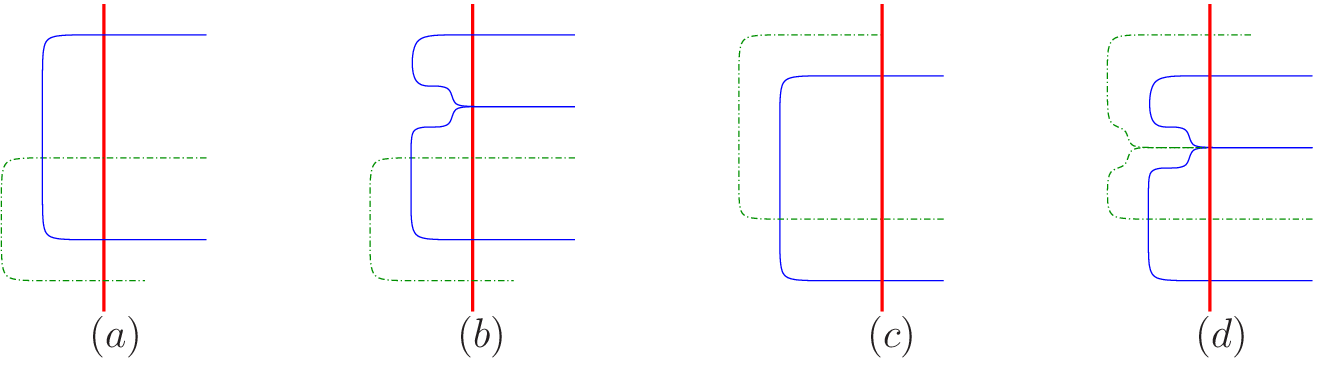}}
\caption{\label{fig:dw_1and2_switched}
{\bf A variation.} On the left of Figures~\ref{fig:dstar_bad} and \ref{fig:dstar_bad_case2}, we might have the picture (c) rather than (a). We then do the finger move in (d) instead of (b). The region on the left is always $D_{z,1}.$}
\end{figure}

Another possible difference is that at the very left of Figures~\ref{fig:dstar_bad}
and \ref{fig:dstar_bad_case2}, the curve $\gamma_{g+1}$ may have an upward rather
than a downward hook, i.e. look as in Figure~\ref{fig:dw_1and2_switched}(c) rather
than (a). If so, instead of the beta finger from the left in Figures~\ref{fig:dstar_bad} and \ref{fig:dstar_bad_case2} (cf. also Figure~\ref{fig:dw_1and2_switched}(b) ), we push a beta-gamma finger as in Figure~\ref{fig:dw_1and2_switched}(d).

Finally, instead of the situations shown in Figures~\ref{fig:dstar_bad} and \ref{fig:dstar_bad_case2}, we might have the same pictures reflected in a horizontal axis. If so, we apply similar finger moves and arrive at the reflections of Figures~\ref{fig:dstar_nice} and \ref{fig:dstar_nice_case2}.

In all cases, the finger moves successfully kill the badness of all regions other than $D_{z,1},$ in which we keep the basepoint $z.$ The result is a nice triple Heegaard diagram for the cobordism $W$.

\subsection{Several two-handle additions} \label {sec:severaltwo}

We now explain how the arguments in this section can be extended to a cobordism $W$ which consists of the addition of several two-handles. We view $W$ as surgery along a link $L \subset Y$ of $l$ components.
 
We start with a multi-pointed Heegaard diagram $(\Sigma, \alphas, \betas, \zees),$ together with another set of basepoints $\wees= \{w_1, \dots, w_l\}$ describing the pair $L \subset Y,$ as in \cite{Links}. Each of the two sets of curves ($\alphas$ and $\betas$) has $g+l-1$ elements.

Applying the algorithm in \cite{SarkarWang} we can make this diagram nice, i.e. such that all regions not containing one of the $z$'s are either bigons or rectangles. For $i=1, \dots, l,$ we denote by $D_{z_i}$ the region containing $z_i.$

\begin{figure}[htbp]
 \center{\includegraphics[width=400pt]{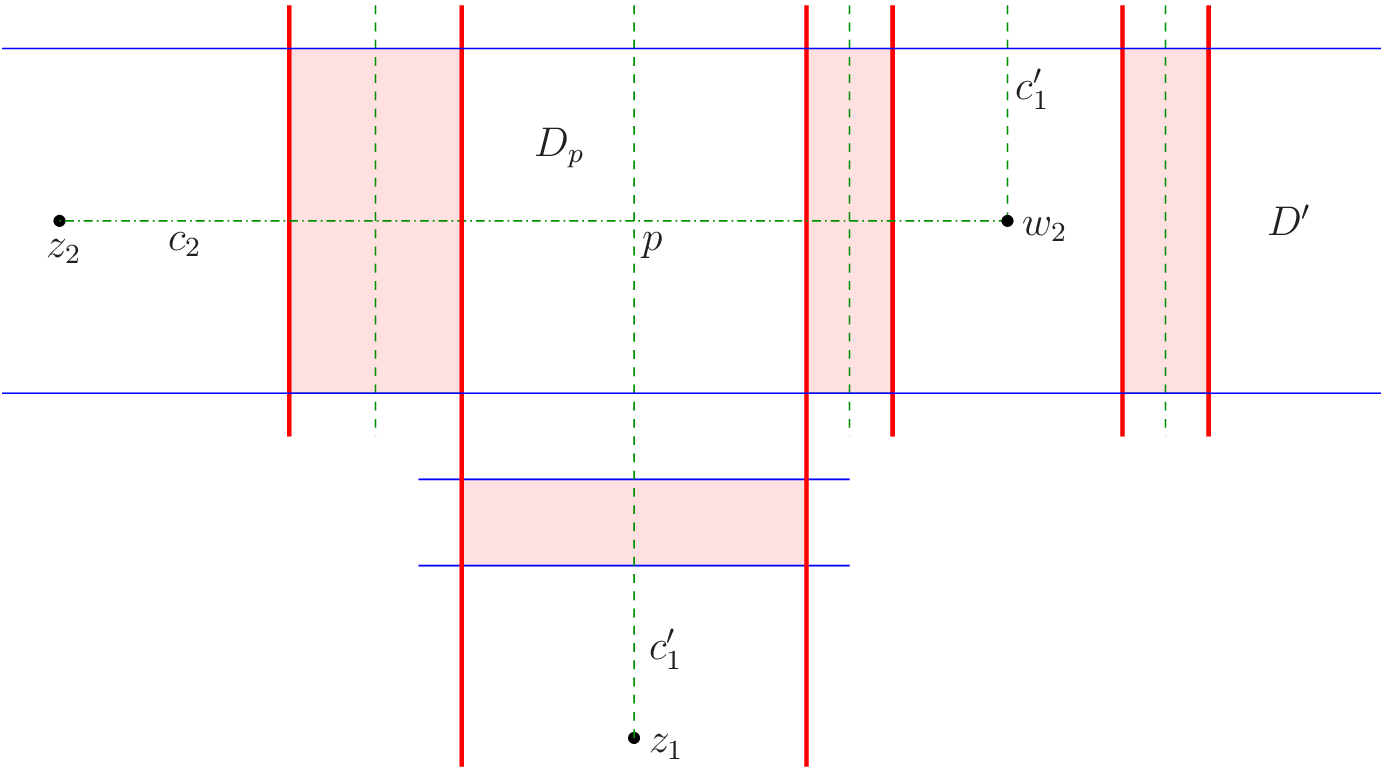}}
\caption{\label{fig:link_embedded_before}
{{\bf Projection of the link to the Heegaard surface.} The point $p$ is the point on the dashed $c'_1$ curve closest to $z_1$ where an intersection with any of the $c_i$ curves (in our case $c_2$) takes place.}} 
\end{figure}

As in Step 1 of Section~\ref{sec:boss}, we inductively remove intersection points between the various components of the projection of $L$ to $\Sigma$. This projection consists of arcs $c_i$ and $c_i'$ with endpoints at $z_i$ and $w_i \ (i=1, \dots, l)$, such that each $c_i$ is disjoint from the beta curves, and each $c_i'$ is disjoint from the alpha curves. Instead of Figure~\ref{fig:knot_embedded_before} we have the situation in Figure~\ref{fig:link_embedded_before}. Again, we stabilize and perform an isotopy to obtain a good diagram with one fewer intersection point, as in Figure~\ref{fig:knot_embedded_after}. Iterating this process (on all link components), we can assume that the projection of $L$ is embedded in the Heegaard surface.

We then add twin gamma curves as in Step 2 of Section~\ref{sec:boss}. For this we need to do several isotopies of the beta curves as in Figure~\ref{fig:adding_gamma_curve_after}. In that figure, if the region on the top left is $D_{z_i},$ the one on the right might be $D_{z_j}$ for $j \neq i$; however, the isotopy can be done as before.

Next, we stabilize the Heegaard diagram $l$ times (once for each link component) to
obtain a triple diagram for the cobordism, as in Step 3 of Section~\ref{sec:boss}. We
then do the analogue of Step 4 by pushing $l$ fingers to kill the badness of the
regions of type $D'$. Since the fingers only pass through rectangles, they do not
intersect each other. The only change is that in Figures~\ref{fig:dprime_ds_bad} and
\ref{fig:dprime_ds_nice}, the region on the very right may contain a different
puncture $z_i$ than the one on the left. By contrast, in Case 2 of Step 4, the region
on the right in Figure~\ref{fig:dprime_bigon} contains the same puncture as the
region on the left, since they lie in the same connected component of the complement
of the beta circles.

At the end of Step 4, the beta curves  split the Heegaard surface into $l$ connected components $C_1, \dots, C_l$. We then do the analogue of Step 5 in Section~\ref{sec:boss}. Note that this step (except for the very last bit, Figure~\ref{fig:dw2_special_bad_region}) only involves moving gamma curves through alpha curves. (Here, we think of the move in Figure~\ref{fig:dw2_special_handleslide} as a single step, rather than as a handleslide followed by an isotopy.) Therefore, we can perform the moves in this step once for each connected component of $L$, independently of each other, because the moves take place in the corresponding component $C_i.$ In the situation considered in Figure~\ref{fig:dw2_special_bad_region}, the gamma curves cross a beta curve; however, the special region $D_s$ is part of a unique $C_i$, so we can perform the isotopy of the gamma curves as before, without interference from another $C_j.$

Finally, for Step 6, note that in all the previous steps we have not destroyed the property that the projection of $L$ to the Heegaard surface is embedded. More precisely, in the part of the stabilized Heegaard diagram shown in Figure~\ref{fig:stabilization}, we take a component of the link projection to be a loop starting in $D_*$, near the upper foot of the handle, going down along $\alpha_{g+1}$ until it reaches $D_{z,1}$, then going inside $D_{z,1}$ until it reaches the intersection of $\gamma_{g+1}$ and $\beta_{g+1},$ and then going along a sub-arc of $\gamma_{g+1}$ to its original departure. These paths remain embedded, and disjoint from each other, throughout Steps 4 and 5. (Indeed, neither $\beta_{g+1}$ nor this sub-arc of $\gamma_{g+1}$ is moved during these steps.) It then suffices to note the finger moves in Step 6 take place in a neighborhood of the projection of the corresponding link component (the path considered above). Therefore, these finger moves can be done without interfering with each other. The result is a nice multi-pointed triple Heegaard diagram for $W$.

\newpage
\section{Computing maps induced by cobordisms}
\label{sec:Last}

Let $W$ be a $4$-dimensional cobordism from $Y_0$ to $Y_3$, and $\ttt$ a
$\spc$-structure on $W$. Choose a self-indexing Morse function on $W$. This
decomposes $W$ as a collection of one-handle additions which together form a
cobordism $W_1$, followed by some two-handle additions forming a cobordism $W_2$, and
three-handle additions forming a cobordism $W_3$, in this order. Let $Y_1$ and $Y_2$
be the intermediate three-manifolds, so that
$$ W = W_1 \ \cup_{Y_1} \ W_2 \ \cup_{Y_2} \ W_3.$$
The map $\Fmap_{W,\ttt}:
\hf(Y_0,\ttt|_{Y_0})\to\hf(Y_3,\ttt|_{Y_3})$
from~\cite{HolDiskFour} is defined as the composition
$\Fmap_{W_3,\ttt|_{W_3}}\circ \Fmap_{W_2,\ttt|_{W_2}}\circ
\Fmap_{W_1,\ttt|_{W_1}}$ of maps associated to each of the pieces $W_1$, $W_2$
and $W_3$. We will review the definitions of these three maps in
Sections~\ref{Section:OneThreeHandles} and~\ref{sec:TwoHandleMaps}. First we
review a few facts about the $H_1(Y; \zz)/\torsion$-action on the hat Heegaard Floer
invariants.

\subsection{$H_1(Y; \zz)/\torsion$-actions} \label{sec:action} 

In \cite[Section 4.2.5]{HolDisk}, Ozsv\'ath and Szab\'o constructed an
action of the group $H_1(Y; \zz)/\torsion \cong \Hom (H^1(Y; \zz),
\zz)$ on $\hf(Y)$. In \cite{HolDiskFour} (see also \cite[Section
2]{OSabs}), they also showed that the cobordism maps $\Fmap_{W, \ttt}$
extend to maps $\Lambda^*(H^1(W;\zz)/\torsion)\otimes\hf(Y_0,\ttt|_{Y_0})\to\hf(Y_3,\ttt|_{Y_3}).$ Moreover, if $W$ is a
cobordism from $Y_1$ to $Y_2$ (endowed with a $\spc$ structure $\ttt$)
and we denote by $j_i: H_1(Y_i; \zz)/\torsion \to H_1(W; \zz)/\torsion \
\ (i=1,2)$ the natural inclusions, then for any $\zeta \in H_1(W;
\zz)/\torsion$ of the form $\zeta = j_1(\zeta_1) - j_2(\zeta_2),  \
\zeta_i \in H_1(Y_i; \zz)/\torsion \ (i=1,2)$ one has
$$ (\zeta \otimes \Fmap_{W, \ttt})(x) = \Fmap_{W, \ttt}(\zeta_1 \cdot x) - \zeta_2 \cdot \Fmap_{W, \ttt}(x).$$

This equality has the following immediate corollaries:
\begin {lemma}
\label {l1}
If $\zeta_1 \in \ker (j_1),$ then $\zeta_1 \cdot x \in \ker ( \Fmap_{W, \ttt})$ for any $x \in \hf (Y_1, \ttt|_{Y_1}).$ 
\end {lemma}

\begin {lemma}
\label {l2}
If $\zeta_2 \in \ker (j_2),$ then $\zeta_2 \cdot y = 0$ for any $y \in \im ( \Fmap_{W, \ttt}).$ 
\end {lemma}

Consider now a 3-manifold of the form $Y\#\left(\#^nS^1\times S^2\right)$, and let $\s_0$ denote the torsion $\spc$-structure on $\#^n S^1\times S^2$. Then for any $\spc$-structure $\s$ on $Y$ there is an isomorphism
\begin{equation}\label{ConSumIsoEq}
\hf(Y,\s)\otimes H_*(T^n)\stackrel{\cong}{\longrightarrow}
\hf\left(Y\#\left(\#^nS^1\times S^2\right),\s\#\s_0\right),
\end{equation} 
as $\F\left[(H_1(Y; \zz)/\torsion)\oplus H_1(\#^nS^1\times S^2; \zz)\right]$-modules, cf. \cite[Theorem  1.5]{HolDiskTwo}. Here, the action of $H_1(\#^nS^1\times S^2; \zz)$ on $\hf(Y,\s)$ is trivial, as is the action of $H_1(Y; \zz)/\torsion$ on $\hf\left(\#^nS^1\times S^2,\s_0\right)$. Further, the action of $H_1(\#^nS^1\times S^2; \zz)/\torsion\cong H^1(T^n; \zz)$ on $\hf(\#^n S^1\times S^2,\s_0)\cong H_*(T^n;\F)$ is exactly the given by cap product $\cap:H^1(T^n;\zz)\otimes H_*(T^n;\F)\to H_{*-1}(T^n;\F)$.

\subsection{Maps associated to one- and three-handle additions}
\label{Section:OneThreeHandles}

Next, we review the definition of the Heegaard Floer maps induced by one- and three-handle additions, cf. \cite[Section 4.3]{HolDiskFour}. 

Suppose that $W_1$ is a cobordism from $Y_0$ to $Y_1$ built entirely from $1$-handles. Let $\ttt$ be a $\spc$-structure on $W_1$. The map $\Fmap_{W_1,\ttt}:\hf(Y_0,\ttt|_{Y_0})\to\hf(Y_1,\ttt|_{Y_1})$ is constructed as follows. If $h_1,h_2, \dots, h_n$ are the $1$-handles in the cobordism, for each $i=1, \dots, n$ pick a path $\xi_i$ in $Y_0$, joining the two feet of the handle $h_i.$ This induces a connected sum decomposition $Y_1\cong Y_0\#(\#^nS^1\times S^2)$, where the first homology of each $S^1 \times S^2$ factor is generated by the union of $\xi_i$ with the core of the corresponding handle. Further, the restriction of $\ttt$ to the $(S^1\times S^2)$-summands in $Y_1$ is torsion. It follows that $\hf(Y_1,\ttt|_{Y_1})\cong \hf(Y_0,\ttt|_{Y_0})\otimes H_*(T^n )$. Let $\theta$ be the generator of the top-graded part of $H_*(T^n)$. Then the Heegaard Floer map induced by $W_1$ is given by
$$\Fmap_{W_1,\ttt}(x)=x\otimes\theta.$$
It is proved in \cite[Lemma 4.13]{HolDiskFour} that, up to composition with canonical isomorphisms, $\Fmap_{W_1,\ttt}$ does not depend on the choices made in its construction, such as the choice of the paths $\xi_i.$

Dually, suppose that $W_3$ is a cobordism from $Y_2$ to $Y_3$ built
entirely from $3$-handles. Let $\ttt$ be a $\spc$-structure on
$W_3$. The map
$\Fmap_{W_3,\ttt}:\hf(Y_2,\ttt|_{Y_2})\to\hf(Y_3,\ttt|_{Y_3})$ is
constructed as follows. One can reverse $W_3$ and view it as attaching
$1$-handles on $Y_3$ to get $Y_2.$ After choosing paths between the
feet of these $1$-handles in $Y_3,$ we obtain a decomposition
$Y_2\cong Y_3\#(\#^mS^1\times S^2)$ (where $m$ is the number of
$3$-handles of $W_3$). Further, the restriction of $\ttt$ to the
$(S^1\times S^2)$-summands in $Y_2$ is torsion. It follows that
$\hf(Y_2,\ttt|_{Y_2})\cong \hf(Y_3,\ttt|_{Y_3})\otimes H_*(T^m)$. Let
$\eta$ be the generator of the lowest-graded part of $H_*(T^m)$. Then
the Heegaard Floer map induced by $W_3$ is given
by $$\Fmap_{W_3,\ttt}(x\otimes\eta)=x$$
and $$\Fmap_{W_3,\ttt}(x\otimes\omega)=0$$ for any homogeneous
generator $\omega$ of $H_*(T^m)$ not lying in the minimal
degree. Again, the map is independent of the choices made in its
construction.

\subsection{Maps associated to two-handle additions}
\label{sec:TwoHandleMaps}
Let $W_2$ be a two-handle cobordism from $Y_1$ to $Y_2$, corresponding
to surgery on a framed link $L$ in $Y_1$, and let $\ttt$ be a
$\spc$-structure on $W_2$. Let $(\Sigma',\alphas',\betas',\gammas',z')$ be
a triple Heegaard diagram subordinate to a bouquet $B(L)$ for $L$, as
in the beginning of Section~\ref{sec:Two}. Then, in particular,
$(\Sigma',\alphas',\betas',z')$ is a Heegaard diagram for $Y_1$,
$(\Sigma',\alphas',\gammas',z')$ is a Heegaard diagram for $Y_2$, and
$(\Sigma',\betas',\gammas',z')$ is a Heegaard diagram for $\#^{g-l}(S^1\times
S^2)$. (Here, $g$ is the genus of $\Sigma'$ and $l$ the number of components of $L$.) The $\spc$-structure $\ttt$ induces a $\spc$-structure
(still denoted $\ttt$) on the four-manifold $W_{\alpha',\beta',\gamma'}$
specified by $(\Sigma',\alphas',\betas',\gammas',z').$ (Note that
$W_{\alpha',\beta',\gamma'}$ can be viewed as a subset of $W_2$.)
Consequently, there is an induced map
\[
\Fmap_{\Sigma',\alphas',\betas',\gammas',z',\ttt}:\hf(Y_1,\ttt|_{Y_1})\otimes\hf\left(\#^{g-l}(S^1\times
S^2),\ttt|_{\#^{g-l}(S^1\times S^2)}\right)\to \hf(Y_2,\ttt|_{Y_2}),
\]
as discussed in Section~\ref{sec:prels}.

The $\spc$-structure $\ttt|_{\#^{g-l}(S^1\times S^2)}$ is necessarily torsion, so
\[\hf\left(\#^{g-l}(S^1\times S^2),\ttt|_{\#^{g-l}(S^1\times S^2)}\right)\cong H_*(T^{g-l})=H_*(S^1)^{\otimes
  (g-l)}.\]
Let $\theta$ denote the generator for the top-dimensional part of
$\hf\left(\#^{g-l}(S^1\times S^2),\ttt|_{\#^{g-l}(S^1\times S^2)}\right)$. Then we define
(cf. \cite[Section 4.1]{HolDiskFour}) the map
\[
\Fmap_{W_2,\ttt}: \hf(Y_1,\ttt|_{Y_1})\to \hf(Y_2,\ttt|_{Y_2})
\]
by $\Fmap_{W_2,\ttt}(x)=\Fmap_{\Sigma',\alphas',\betas',\gammas',z',\ttt}(x\otimes\theta)$.

Now, consider instead the nice, $l$-pointed triple Heegaard diagram
$(\Sigma,\alphas,\betas,\gammas,\zees)$ constructed in
Section~\ref{sec:Two}. As discussed in the beginning of
Section~\ref{sec:Two}, $(\Sigma,\alphas,\betas,\gammas,\zees)$ is
strongly equivalent to a split triple Heegaard diagram whose reduction
$(\Sigma',\alphas',\betas',\gammas',z')$ is subordinate to a bouquet
$B(L)$ as above. Let $\Theta$ be the generator for the top-dimensional
part of $\hf(\Sigma,\betas,\gammas,\zees)\cong H_*(T^{g+l-1})$. Then,
by Diagram~\eqref{diagram:commutes}, with $k=l-1$, we have
\[
\rank\left(\Fmap_{W_2,\ttt}\right)=\frac{1}{2^{l-1}}\rank\left(\Fmap_{\Sigma,\alphas,\betas,\gammas,\zees,\ttt}(\cdot\otimes\Theta)\right).
\]

In light of Corollary~\ref{lemma:admissible} and
Proposition~\ref{EasternOrthodox}, the rank of the map
$\Fmap_{\Sigma,\alphas,\betas,\gammas,\zees,\ttt}$ can be computed
combinatorially. Further, since the triple Heegaard diagram
$(\Sigma,\alphas,\betas,\gammas,\zees)$ is nice, so are each of the
three (ordinary) Heegaard diagrams it specifies. Consequently,
by~\cite{SarkarWang}, the element
$\Theta\in\hf(\Sigma,\betas,\gammas,\zees)$ can be explicitly
identified (as can a representative for $\Theta$ in
$\widehat{CF}(\Sigma,\betas,\gammas,\zees)$). Therefore, the rank of
$\Fmap_{W_2,\ttt}$ can be computed combinatorially.

\subsection{Putting it all together}

Recall that $W$ is a $4$-dimensional cobordism from $Y_0$ to $Y_3$, $\ttt$ a
$\spc$-structure on $W$, and that $W$ is decomposed as a collection of one-handle
additions $W_1$, followed by some two-handle additions $W_2$, and three-handle
additions $W_3$, with $Y_1$ and $Y_2$ the intermediate three-manifolds, so that 
$$ W = W_1 \ \cup_{Y_1} \ W_2 \ \cup_{Y_2} \ W_3.$$

As in Section~\ref{sec:action}, we consider the maps 
\begin {equation}
\label {jays}
 j_1: H_1(Y_0; \zz)/\torsion \to H_1(W; \zz)/\torsion, \ \ \  j_2: H_1(Y_3; \zz)/\torsion \to H_1(W; \zz)/\torsion.
\end {equation}

\begin {lemma}
\label {onto1}
If $j_1$ is surjective, then $\im (\Fmap_{W_2, \ttt|_{W_2}} \circ \Fmap_{W_1, \ttt|_{W_1}}) = \im (\Fmap_{W_2, \ttt|_{W_2}}).$ 
\end {lemma}

\begin {proof}
The cobordism $W_1$ consists of the addition of some $1$-handles $h_1,
\dots, h_n.$ As in Section~\ref{Section:OneThreeHandles}, we choose
paths $\xi_i$ in $Y_0$ joining the two feet of the handle $h_i.$ The
union of $\xi_i$ with the core of $h_i$ produces a curve in $W_1$,
which in turn gives an element $e_i \in H_1(W_1 \cup W_2;
\zz)/\torsion.$ Since $W$ is obtained from $W_1\cup W_2$ by adding
$3$-handles, we have $H_1(W; \zz)\cong H_1(W_1\cup W_2; \zz)$, so the
hypothesis implies that the map
$$j_1' : H_1(Y_0; \zz)/\torsion \to H_1(W_1 \cup W_2; \zz)/\torsion,$$
is surjective. Hence there exist disjoint, embedded curves $c_i$ in
$Y_0$ (disjoint from all the $\xi_j$) such that $j'_1([c_i]) = -e_i,
i=1, \dots, n.$ We can connect sum $\xi_i$ and $c_i$ to get new paths
$\xi_i'$ in $Y_0$ between the two feet of $h_i.$ Using the paths
$\xi_i'$ we get a connected sum decomposition $Y_1 \cong
Y_0\#(\#^nS^1\times S^2)$ as in Section~\ref{Section:OneThreeHandles},
with the property that the inclusion of the summand $H_1(\#^nS^1\times
S^2; \zz) \subset H_1(Y_1; \zz)$ in $H_1(W_1 \cup W_2; \zz)/\torsion$
is trivial. Since we can view $W_1\cup W_2$ as obtained from $W_2$ by
adding $3$-handles (which do not affect $H_1$), it follows that the inclusion of  $H_1(\#^nS^1\times S^2; \zz)$ in $H_1(W_2; \zz)/\torsion$ is trivial. Lemma~\ref{l1} then says that
$$ \Fmap_{W_2, \ttt|_{W_2}} (\zeta \cdot x) = 0,$$ 
for any $\zeta \in H_1(\#^nS^1\times S^2; \zz), x \in \hf(Y_1, \ttt|_{Y_1}) \cong \hf(Y_0, \ttt|_{Y_0}) \otimes H_*(T^n).$ Thus the kernel of $\Fmap_{W_2, \ttt|_{W_2}}$ contains all elements of the form $y \otimes \omega$, where $y \in \hf(Y_0, \ttt|_{Y_0})$ and $\omega \in H_*(T^n)$ is any homogeneous element not lying in the top grading of $H_*(T^n ).$ On the other hand, from Section~\ref{Section:OneThreeHandles} we know that the image of $\Fmap_{W_1, \ttt|_{W_1}}$ consists exactly of the elements $y \otimes \theta,$ where $\theta$ is the top degree generator of $H_*(T^n).$ Therefore,
$$ \im (\Fmap_{W_1, \ttt|_{W_1}}) + \ker (\Fmap_{W_2, \ttt|_{W_2}}) = \hf(Y_1, \ttt|_{Y_1}).$$
This gives the desired result.
\end {proof}

\begin {lemma}
\label {onto2}
If $j_2$ is surjective, then $\im (\Fmap_{W_2, \ttt|_{W_2}}) \cap \ker (\Fmap_{W_3, \ttt|_{W_3}}) = 0.$ 
\end {lemma}

\begin {proof} This is similar to the proof of Lemma~\ref{onto1}. A suitable choice of paths enables us to view $Y_2$ as $Y_3 \#(\#^mS^1\times S^2),$ such that the inclusion of the summand $H_1(\#^mS^1\times S^2; \zz) \subset H_1(Y_2; \zz)$ in $H_1(W_2; \zz)/\torsion$ is trivial. Lemma~\ref{l2} then says that $\zeta \cdot y = 0$ for any $y \in \im ( \Fmap_{W_2, \ttt|_{W_2}})$ and $\zeta \in  H_1(\#^mS^1\times S^2).$ In other words, every element in the image of $\Fmap_{W_2, \ttt|_{W_2}}$ must be of the form $y = x \otimes \eta,$ where $x \in \hf(Y_3, \ttt|_{Y_3})$ and $\eta$ is the lowest degree generator of $H_*(T^n).$ On the other hand, from Section~\ref{Section:OneThreeHandles} we know that the kernel of the map $\Fmap_{W_3, \ttt|_{W_3}}$ does not contain any nonzero elements of the form $x \otimes \eta.$ 
\end {proof}

\begin{theorem}Let $W$ be a cobordism from $Y_0$ to $Y_3$, and $\ttt$ a $\spc$-structure on $W$. Assume that the maps $j_1$ and $j_2$ from Formula~\eqref{jays} are surjective. Then in each (relative) grading $i$ the rank of $\Fmap_{W,\ttt}:\hf_i(Y_0,\ttt|_{Y_0})\to\hf_*(Y_3,\ttt|_{Y_3})$ can be computed combinatorially.
\end{theorem}

\begin{proof}
The map $\Fmap_{W,\ttt}$ is, by definition, the composition $\Fmap_{W_3,\ttt|_{W_3}}\circ \Fmap_{W_2,\ttt|_{W_2}}\circ \Fmap_{W_1,\ttt|_{W_1}}$. Lemmas~\ref{onto1} and \ref{onto2} imply that
$$ \im (\Fmap_{W_2, \ttt|_{W_2}} \circ  \Fmap_{W_1, \ttt|_{W_1}} ) \cap \ker (\Fmap_{W_3, \ttt|_{W_3}}) = 0$$
or, equivalently, 
$$ \rk (\Fmap_{W_3, \ttt|_{W_3}} \circ \Fmap_{W_2, \ttt|_{W_2}} \circ  \Fmap_{W_1, \ttt|_{W_1}} ) = \rk (\Fmap_{W_2, \ttt|_{W_2}} \circ  \Fmap_{W_1, \ttt|_{W_1}} ).$$

Using Lemma~\ref{onto1} again, the expression on the right is the same as the rank of $\Fmap_{W_2, \ttt|_{W_2}}.$
Thus, the maps $\Fmap_{W,\ttt}$ and $\Fmap_{W_2, \ttt|_{W_2}}$ have the same rank. As explained in Section~\ref{sec:TwoHandleMaps}, the rank of  $\Fmap_{W_2, \ttt|_{W_2}}$ can be computed combinatorially. Note that the relative gradings on the generators of the chain complexes are also combinatorial, using the formula for the Maslov index in \cite[Corollary 4.3]{Lipshitz}. This completes the proof.
\end{proof}

\begin {remark} 
In fact, using Sarkar's remarkable formula for the Maslov index of triangles
\cite[Theorem 4.1]{SarkarIndex}, the absolute gradings on the Heegaard Floer complexes
can be computed combinatorially, and so the rank of $\Fmap_{W,\ttt}$ in each absolute grading can be computed as well.
\end {remark}

\newpage
\section{An example}
\label{Sec:Example}

\begin{figure}[htbp]
 \center{\includegraphics[width=250pt]{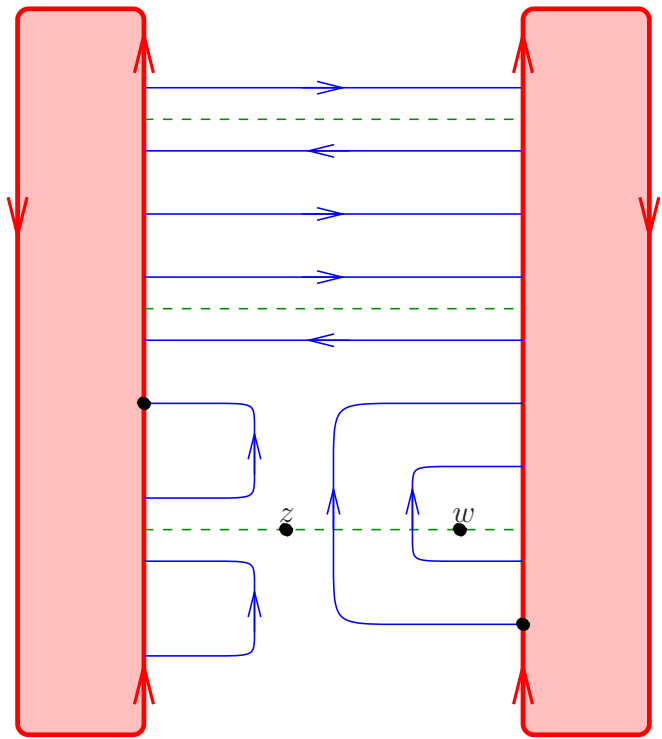}} 
\caption{\label{fig:poincare_trefoil}
 {\bf A nice Heegaard diagram for the trefoil knot.} The thick curves are 
alpha curves and the thin ones are beta curves. The two shaded areas in 
the diagram are glued together via a reflection and rotation, such that 
the two black dots on the alpha curves are identified. The knot is given 
by the dashed curves. Its projection is already embedded in the Heegaard 
surface.} 
\end {figure}

\begin{figure}[htbp]
 \center{\includegraphics[width=440pt]{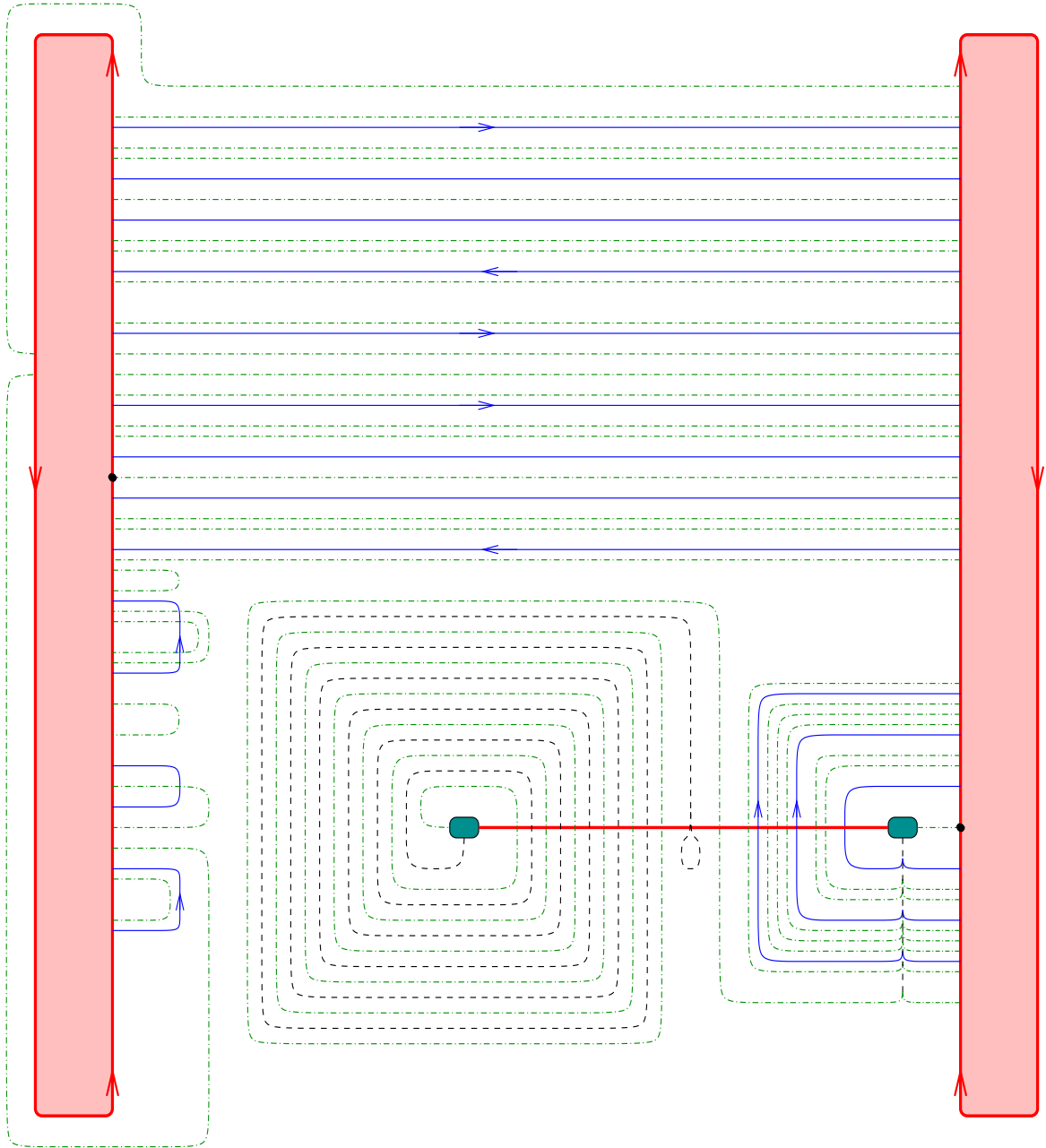}} 
\caption{\label{fig:poincare_nice} 
{\bf A nice triple Heegaard diagram for $+1$ surgery on the right-handed 
trefoil.} The thick curves are alpha curves, the thin curves are beta 
curves, and the interrupted curves are gamma curves. The lightly shaded 
areas are identified as in Figure~\ref{fig:poincare_trefoil}, while the 
darkly shaded areas represent the new handle which has been attached. 
Again, we are using the train-track convention as in Figure 
\ref{fig:train-track}. The basepoint is the point at infinity.} 
\end 
{figure}

We give a nice triple Heegaard diagram for the cobordism from the 
three-sphere to the Poincar\'e homology sphere, viewed as the $+1$ 
surgery on the right-handed trefoil. The right-handed trefoil knot admits 
the nice Heegaard diagram shown in Figure~\ref{fig:poincare_trefoil}, 
which is isotopic to \cite[Figure 14]{SarkarWang}. Applying the algorithm 
described in Section~\ref{sec:Two}, we obtain the nice triple Heegaard 
diagram shown in Figure~\ref{fig:poincare_nice}. We leave the actual 
computation of the cobordism map to the interested reader.

\end{document}